\documentclass[12pt,seceqn,secthm]{elsart}%
\usepackage{amsfonts}
\usepackage{amsmath}
\usepackage{amssymb}
\usepackage[dvips]{graphicx}%
\renewcommand{\theenumi}{\roman{enumi}}

\makeatletter\def\Tiny{\@setfontsize\Tiny\@vpt{6}}\makeatother
\newcommand{\absolutelycontinuous}{<\mathrel{\vcenter{\llap{\Tiny
$\boldsymbol{<}\mkern0.8mu\joinrel$}}}\mathrel{\vcenter{\llap{\Tiny
\raisebox{0.05pt}{\llap{\Tiny
$\boldsymbol{<}\mkern0.6mu\joinrel$}}\raisebox{-0.05pt}{\llap{\Tiny
$\boldsymbol{<}\mkern0.6mu\joinrel$}}}}}\mathrel{\vcenter{\llap{\Tiny
$\boldsymbol{<}\mkern0.4mu\joinrel$}}}}
\begin{document}
%
\begin{frontmatter}%
%

\title{Iterated function systems\protect\rule{0pt}{96pt}%
, representations, and Hilbert space}%
%

\author{Palle E. T. Jorgensen\thanksref{label1}}\thanks
[label1]{This material is based upon work supported by the U.S.
National Science Foundation under Grant No.\ DMS-0139473 (FRG).}%

%

\address{Department of Mathematics,
The University of Iowa,
14 MacLean Hall,
Iowa City, IA 52242-1419,
U.S.A.}%
\ead{jorgen@math.uiowa.edu}%
\ead[url]{http://www.math.uiowa.edu/\symbol{126}jorgen}%
%

\begin{abstract}
In this paper, we are concerned with spectral-theoretic features of general
iterated function systems (IFS). Such systems arise from the study of
iteration limits of a finite family of maps $\tau_{i}$, $i=1,\dots,N$, in some
Hausdorff space $Y$. There is a standard construction which generally allows
us to reduce to the case of a compact invariant subset $X\subset Y$.
Typically, some kind of contractivity property for the maps $\tau_{i}$ is
assumed, but our present considerations relax this restriction. This means
that there is then not a natural equilibrium measure $\mu$ available which
allows us to pass the point-maps $\tau_{i}$ to operators on the Hilbert space
$L^{2}\left(  \mu\right)  $. Instead, we show that it is possible to realize
the maps $\tau_{i}$ quite generally in Hilbert spaces $\mathcal{H}\left(
X\right)  $ of square-densities on $X$. The elements in $\mathcal{H}\left(
X\right)  $ are equivalence classes of pairs $\left(  \varphi,\mu\right)  $,
where $\varphi$ is a Borel function on $X$, $\mu$ is a positive Borel measure
on $X$, and $\int_{X}\left\vert \varphi\right\vert ^{2}\,d\mu<\infty$. We say
that $\left(  \varphi,\mu\right)  \sim\left(  \psi,\nu\right)  $ if there is a
positive Borel measure $\lambda$ such that $\mu\absolutelycontinuous\lambda$,
$\nu\absolutelycontinuous\lambda$, and%
\[
\varphi\,\sqrt{\frac{d\mu}{d\lambda}}=\psi\,\sqrt{\frac{d\nu}{d\lambda}%
}\,,\qquad\lambda\;\mathrm{a.e.}\text{ on }X.
\]
We prove that, under general conditions on the system $\left(  X,\tau
_{i}\right)  $, there are isometries%
\[
S_{i}\colon\left(  \varphi,\mu\right)  \longmapsto\left(  \varphi\circ
\sigma,\mu\circ\tau_{i}^{-1}\right)
\]
in $\mathcal{H}\left(  X\right)  $ satisfying $\sum_{i=1}^{N}S_{i}S_{i}^{\ast
}=I={}$the identity operator in $\mathcal{H}\left(  X\right)  $. For the
construction we assume that some mapping $\sigma\colon X\rightarrow X$
satisfies the conditions $\sigma\circ\tau_{i}=\operatorname*{id}\nolimits_{X}%
$, $i=1,\dots,N$.

We further prove that this representation in the Hilbert space $\mathcal{H}%
\left(  X\right)  $ has several universal properties.%
\end{abstract}%
%

\begin{keyword}
Hilbert space, Cuntz algebra, completely positive map, creation operators,
wavelet packets, pyramid algorithm, product measures, orthogonality relations,
equivalence of measures, iterated function systems (IFS), scaling function,
multiresolution, subdivision scheme, singular measures, absolutely continuous
measures
\renewcommand{\MSC}{{\par\leavevmode\hbox{\it2000 MSC:\ }}}\MSC
42C40; 42A16; 43A65; 42A65%
\end{keyword}%
%

\end{frontmatter}%

\section{\label{Int}Introduction}

In this paper we are concerned with \emph{iterated function systems} (IFS) and
their representation in Hilbert space. For contractive IFS's, there is a known
standard construction of a family of measures, and Hilbert spaces induced by
these measures. However, the constructions are not universal in any reasonable
sense, and they only admit a very restricted family of covariant measures.

Let $X$ be a compact metric space, and let $\tau_{i}\colon X\rightarrow X$,
$i=1,\dots,N$, satisfy%
\begin{equation}
d\left(  \tau_{i}\left(  x\right)  ,\tau_{i}\left(  y\right)  \right)  \leq
Cd\left(  x,y\right)  ,\qquad i=1,\dots,N,\;x,y\in X, \label{eqInt.1}%
\end{equation}
for some $C$, $0<C<1$. Let $p_{i}>0$ be given such that $\sum_{i=1}^{N}%
p_{i}=1$. Then it follows from a theorem of Hutchinson \cite{Hut81} that there
is a unique positive Borel measure $\mu=\mu_{\left(  p\right)  }$ on $X$ such
that $\mu\left(  X\right)  =1$, and%
\begin{equation}
\sum_{i=1}^{N}p_{i}\;\mu\circ\tau_{i}^{-1}=\mu, \label{eqInt.2}%
\end{equation}
where the measures $\mu\circ\tau_{i}^{-1}$ are defined by $\mu\circ\tau
_{i}^{-1}\left(  E\right)  :=\mu\left(  \tau_{i}^{-1}\left(  E\right)
\right)  $, $E\in\mathcal{B}\left(  X\right)  ={}$the Borel subsets of $X$,
where%
\begin{equation}
\tau_{i}^{-1}\left(  E\right)  :=\left\{  \,x\in X\mid\tau_{i}\left(
x\right)  \in E\,\right\}  . \label{eqInt.3}%
\end{equation}
We shall need a \textquotedblleft variable-coefficient
version\textquotedblright\ of (\ref{eqInt.2}) which is motivated by
applications to wavelets; see \cite{Jor04a} and \cite{BrJo02b}. In this
version of (\ref{eqInt.2}), there is a whole family of measures $\mu_{f}$
indexed by vectors $f$ in some complex Hilbert space $\mathcal{K}$, and
moreover there is a finite family of isometries $S_{i}\colon\mathcal{K}%
\rightarrow\mathcal{K}$, $i=1,\dots,N$, such that%
\begin{equation}
\sum_{i=1}^{N}S_{i}S_{i}^{\ast}=I_{\mathcal{K}}\,, \label{eqInt.4}%
\end{equation}
and (\ref{eqInt.2}) takes the form%
\begin{equation}
\sum_{i=1}^{N}\mu_{S_{i}^{\ast}f}\circ\tau_{i}^{-1}=\mu_{f}\,. \label{eqInt.5}%
\end{equation}
Isometries $S_{i}$ subject to (\ref{eqInt.4}) are said to satisfy the Cuntz
relations, or to define a representation of the Cuntz algebra $\mathcal{O}%
_{N}$; see \cite{Cun77}. The algebra $\mathcal{O}_{N}$ is a simple $C^{\ast}%
$-algebra, and its representations are ubiquitous in analysis and applied
mathematics. A special class of these relations is known to define subband
filters in signal processing, to define subdivision algorithms in computer
graphics, and to define effective pyramid algorithms in wavelet analysis; see
\cite{BrJo02b}, \cite{Jor03}. However, the classical approach to subdivisions
via (\ref{eqInt.2}) is known not to suffice for the representation of wavelet
systems, see \cite{Jor01b}, not even for the simplest quadrature mirror
filters which are used for the standard Haar wavelet or for the Daubechies wavelets.

Readers not familiar with wavelets may pick up the essentials from Chapter 5
of the classic \cite{Dau92}. More current results, presented from the
viewpoint of operator theory, may be found in Chapter 2 of the monograph
\cite{BrJo02b}, or in the survey paper \cite{Jor03}.

Another aim of the present paper is to relax the contractivity condition
(\ref{eqInt.1}). Our starting point is a compact Hausdorff space $X$ and
continuous maps $\sigma\colon X\rightarrow X$, $\tau_{i}\colon X\rightarrow
X$, $i=1,\dots,N$, such that
\begin{equation}
\sigma\circ\tau_{i}=\operatorname*{id}\nolimits_{X}\,. \label{eqInt.6}%
\end{equation}
It follows from (\ref{eqInt.6}) that $\sigma$ is onto, and that each $\tau
_{i}$ is one-to-one. We will be especially interested in the case when there
are distinct branches $\tau_{i}\colon X\rightarrow X$ such that%
\begin{equation}
\bigcup_{i=1}^{N}\tau_{i}\left(  X\right)  =X. \label{eqInt.7}%
\end{equation}
For such systems, we show in Section \ref{Fro} that there is a
\emph{universal} representation of $\mathcal{O}_{N}$ in a Hilbert space
$\mathcal{H}\left(  X\right)  $ which is functorial, is naturally defined, and
contains every representation of $\mathcal{O}_{N}$.

The elements in the universal Hilbert space $\mathcal{H}\left(  X\right)  $
are equivalence classes of pairs $\left(  \varphi,\mu\right)  $ where
$\varphi$ is a Borel function on $X$ and where $\mu$ is a positive Borel
measure on $X$. We will set $\varphi\sqrt{d\mu}:=\operatorname*{class}\left(
\varphi,\mu\right)  $ for reasons which we spell out below.

While our present methods do adapt to the more general framework when the
space $X$ of (\ref{eqInt.6})--(\ref{eqInt.7}) is not assumed compact, but only
$\sigma$-compact, we will still restrict the discussion here to the compact
case. This is for the sake of simplicity of the technical arguments. But we
encourage the reader to follow our proofs below, and to formulate for
him/her\-self the corresponding results when $X$ is not necessarily assumed
compact. Moreover, if $X$ is not compact, then there is a variety of special
cases to take into consideration, various abstract notions of ``escape to
infinity''. We leave this wider discussion for a later investigation, and we
only note here that our methods allow us to relax the compactness restriction
on $X$.

There is a classical construction in operator theory which lets us realize
point transformations in Hilbert space. It is called the Koopman
representation; see, for example, \cite[p.~135]{Mac89}. But this approach only
applies if the existence of invariant, or quasi-invariant, measures is
assumed. In general such measures are not available. The present paper
proposes a different way of realizing families of point transformations in
Hilbert space in a general context where no such assumptions are made. Our
Hilbert spaces are motivated by a construction due to S.~Kakutani
\cite{Kak48}, L.~Schwartz, and E.~Nelson \cite{Nel69}, among others. The
reader is also referred to an updated presentation of the measure-class
Hilbert spaces due to Tsirelson \cite{Tsi03} and Arveson \cite[Chapter
14]{Arv03a}.

We say that $\left(  \varphi,\mu\right)  \sim\left(  \psi,\nu\right)  $ if
there is a third positive Borel measure $\lambda$ on $X$ such that
$\mu\absolutelycontinuous\lambda$, $\nu\absolutelycontinuous\lambda$, and%
\begin{equation}
\varphi\,\sqrt{\frac{d\mu}{d\lambda}}=\psi\,\sqrt{\frac{d\nu}{d\lambda}%
}\,,\qquad\lambda\;\mathrm{a.e.}\text{ on }X, \label{eqInt.8}%
\end{equation}
where $\absolutelycontinuous$ denotes relative absolute continuity, and where
$d\mu/d\lambda$ denotes the usual Radon-Nikodym derivative, i.e.,
$d\mu/d\lambda\in L^{1}\left(  \lambda\right)  $, and $d\mu=\left(
d\mu/d\lambda\right)  \,d\lambda$.

In Section \ref{Iso}, we review some basic properties of the Hilbert space
$\mathcal{H}\left(  X\right)  $. This space is called the Hilbert space of
$\sigma$-functions, or square densities, and it was studied for different
reasons in earlier papers of L.~Schwartz, E.~Nelson \cite{Nel69}, and
W.~Arveson \cite{Arv03b}.

Our first new result is the fact that the isometries $S_{i}\colon
\mathcal{H}\left(  X\right)  \rightarrow\mathcal{H}\left(  X\right)  $ are
defined by
\begin{equation}
S_{i}\colon\left(  \varphi,\mu\right)  \longmapsto\left(  \varphi\circ
\sigma,\mu\circ\tau_{i}^{-1}\right)  _{\mathstrut}, \label{eqInt.9}%
\end{equation}
or equivalently, $S_{i}\colon\varphi\sqrt{d\mu}\mapsto\varphi\circ
\sigma\;\sqrt{d\mu\circ\tau_{i}^{-1}}$, and that the operators satisfy the
Cuntz relations (\ref{eqInt.4}).

Note that, at the outset, it is not even clear \emph{a priori} that $S_{i}$ in
(\ref{eqInt.9}) defines a transformation of $\mathcal{H}\left(  X\right)  $.
To verify this, we will need to show that if two equivalent pairs are
substituted on the left-hand side in (\ref{eqInt.9}), then they produce
equivalent pairs as output, on the right-hand side. Recalling the definition
(\ref{eqInt.8}) of the equivalence relation $\sim$, there is no obvious or
intuitive reason for why this should be so.

To stress the intrinsic transformation rules of $\mathcal{H}\left(  X\right)
$, the vectors in $\mathcal{H}\left(  X\right)  $ are usually denoted
$\varphi\sqrt{d\mu}$ rather than $\left(  \varphi,\mu\right)  $. This is a
suggestive notation which motivates the definition of the inner product of
$\mathcal{H}\left(  X\right)  $. It is also helpful in understanding Theorem
\ref{ThmFro.2} below. If $\varphi\sqrt{d\mu}$ and $\psi\sqrt{d\nu}$ are in
$\mathcal{H}\left(  X\right)  $, we define their Hilbert inner product by%
\begin{equation}
\left\langle \,\varphi\sqrt{d\mu}\bigm|\psi\sqrt{d\nu}\,\right\rangle
:=\int_{X}\bar{\varphi}\,\psi\,\sqrt{\frac{d\mu}{d\lambda}}~\sqrt{\frac{d\nu
}{d\lambda}}\,d\lambda, \label{eqInt.10}%
\end{equation}
where $\lambda$ is some positive Borel measure, chosen such that
$\mu\absolutelycontinuous\lambda$ and $\nu\absolutelycontinuous\lambda$. For
example, we could take $\lambda=\mu+\nu$. To be in $\mathcal{H}\left(
X\right)  $, $\varphi\sqrt{d\mu}$ must satisfy%
\begin{equation}
\left\Vert \varphi\sqrt{d\mu}\right\Vert ^{2}=\int_{X}\left\vert
\varphi\right\vert ^{2}\,\frac{d\mu}{d\lambda}\,d\lambda=\int_{X}\left\vert
\varphi\right\vert ^{2}\,d\mu<\infty. \label{eqInt.11}%
\end{equation}

\section{\label{Iso}Isometries in $\mathcal{H}\left(  X\right)  $}

In this preliminary section we prove three general facts about the process of
inducing operators in the Hilbert space $\mathcal{H}\left(  X\right)  $ from
underlying point transformations in $X$. The starting point is a given
continuous mapping $\sigma\colon X\rightarrow X$, mapping onto $X$.$\ We$ will
be concerned with the special case when $X$ is a compact Hausdorff space, and
when there is one or more continuous branches $\tau_{i}\colon X\rightarrow X$
of the inverse, i.e., when%
\begin{equation}
\sigma\circ\tau_{i}=\operatorname*{id}\nolimits_{X}. \label{eqIsoNew.1}%
\end{equation}
Recall that elements in $\mathcal{H}\left(  X\right)  $ are equivalence
classes of pairs $\left(  \varphi,\mu\right)  $ where $\varphi$ is a Borel
function on $X$, $\mu$ is a positive Borel measure on $X$, and $\int
_{X}\left\vert \varphi\right\vert ^{2}\,d\mu<\infty$. An equivalence class
will be denoted $\varphi\sqrt{d\mu}$, and we show that there are isometries
\begin{equation}
S_{i}\colon\varphi\sqrt{d\mu}\longmapsto\varphi\circ\sigma\;\sqrt{d\mu
\circ\tau_{i}^{-1}}, \label{eqIsoNew.2}%
\end{equation}
with orthogonal ranges in the Hilbert space $\mathcal{H}\left(  X\right)  $.
Moreover, we calculate an explicit formula for the adjoint co-isometries
$S_{i}^{\ast}$.

In the next section, we shall then restrict the setting to the special case of
measures $\mu$ such that $\mu\circ\tau_{i}^{-1}\absolutelycontinuous\mu$,
where $\absolutelycontinuous$ stands for \textquotedblleft absolutely
continuous with respect to\textquotedblright.

\begin{lem}
\label{LemIso.2}Let $X$ be a compact Hausdorff space, and let the mapping
$\sigma\colon X\rightarrow X$ be onto. Suppose $\tau\colon X\rightarrow X$
satisfies $\sigma\circ\tau=\operatorname*{id}\nolimits_{X}$. Assume that both
$\sigma$ and $\tau$ are continuous. Let $\mathcal{H}=\mathcal{H}\left(
X\right)  $ be the Hilbert space of classes $\left(  \varphi,\mu\right)  $
where $\varphi$ is a Borel function on $X$ and $\mu$ is a positive Borel
measure such that $\int\left\vert \varphi\right\vert ^{2}\,d\mu<\infty$. The
equivalence relation is defined in the usual way: two pairs $\left(
\varphi,\mu\right)  $ and $\left(  \psi,\nu\right)  $ are said to be
equivalent, written $\left(  \varphi,\mu\right)  \sim\left(  \psi,\nu\right)
$, if for some positive measure $\lambda$, $\mu\absolutelycontinuous\lambda$,
$\nu\absolutelycontinuous\lambda$, we have the following identity:%
\begin{equation}
\varphi\,\sqrt{\frac{d\mu}{d\lambda}}=\psi\,\sqrt{\frac{d\nu}{d\lambda}}%
\qquad(\mathrm{a.e.}\;\lambda). \label{eqIso.1}%
\end{equation}
Then there is an isometry $S\colon\mathcal{H}\rightarrow\mathcal{H}$ which is
well defined by the assignment%
\begin{equation}
S\left(  \left(  \varphi,\mu\right)  \right)  :=\left(  \varphi\circ\sigma
,\mu\circ\tau^{-1}\right)  , \label{eqIso.2}%
\end{equation}
or%
\[
S\colon\varphi\sqrt{d\mu}\longmapsto\varphi\circ\sigma\;\sqrt{d\mu\circ
\tau^{-1}},
\]
where $\mu\circ\tau^{-1}\left(  E\right)  :=\mu\left(  \tau^{-1}\left(
E\right)  \right)  $, and $\tau^{-1}\left(  E\right)  :=\left\{  \,x\in
X\mid\tau\left(  x\right)  \in E\,\right\}  $, for $E\in\mathcal{B}\left(
X\right)  $.
\end{lem}

\begin{pf}
We leave the verification of the following four facts to the reader; see also
\cite{Nel69}.

\begin{enumerate}
\item \label{LemIso.2proof(1)}If $\displaystyle\varphi\,\sqrt{\frac{d\mu
}{d\lambda}}=\psi\,\sqrt{\frac{d\nu}{d\lambda}}$ for some $\lambda$ such that
$\mu\absolutelycontinuous\lambda$, $\nu\absolutelycontinuous\lambda$, and if
some other measure $\lambda^{\prime}$ satisfies $\mu
\absolutelycontinuous\lambda^{\prime}$, $\nu\absolutelycontinuous\lambda
^{\prime}$, then
\[
\varphi\,\sqrt{\frac{d\mu}{d\lambda^{\prime}}}=\psi\,\sqrt{\frac{d\nu
}{d\lambda^{\prime}}}\qquad(\mathrm{a.e.}\,\lambda^{\prime}).
\]

\item \label{LemIso.2proof(2)}The \textquotedblleft vectors\textquotedblright%
\ in $\mathcal{H}$ are equivalence classes of pairs $\left(  \varphi
,\mu\right)  $ as described in the statement of the lemma. For two elements
$\left(  \varphi,\mu\right)  $ and $\left(  \psi,\nu\right)  $ in
$\mathcal{H}$, define the sum by%
\begin{equation}
\left(  \varphi,\mu\right)  +\left(  \psi,\nu\right)  :=\left(  \phi
\,\sqrt{\frac{d\mu}{d\lambda}}+\psi\,\sqrt{\frac{d\nu}{d\lambda}}%
,\;\lambda\right)  , \label{eqIso.3}%
\end{equation}
where $\lambda$ is a positive Borel measure satisfying $\mu
\absolutelycontinuous\lambda$, $\nu\absolutelycontinuous\lambda$. The sum in
(\ref{eqIso.3}) is also written $\varphi\sqrt{d\mu}+\psi\sqrt{d\nu}$. The
definition of the sum (\ref{eqIso.3}) passes through the equivalence relation
$\sim$, i.e., we get an equivalent result on the right-hand side in
(\ref{eqIso.3}) if equivalent pairs are used as input on the left-hand side. A
similar conclusion holds for the definition (\ref{eqIso.4}) below of the inner
product $\left\langle \,\cdot\mid\cdot\,\right\rangle $ in the Hilbert space
$\mathcal{H}$.

\item \label{LemIso.2proof(3)}Scalar multiplication, $c\in\mathbb{C}$, is
defined by $c\left(  \varphi,\mu\right)  :=\left(  c\mkern2mu\varphi
,\mu\right)  $, and the Hilbert space inner product is
\begin{equation}
\left\langle \,\varphi\sqrt{d\mu}\Bigm|\psi\sqrt{d\nu}\,\right\rangle
=\left\langle \,\left(  \varphi,\mu\right)  \mid\left(  \psi,\nu\right)
\,\right\rangle :=\int_{X}\bar{\varphi}\,\psi\,\sqrt{\frac{d\mu}{d\lambda}%
}\,\sqrt{\frac{d\nu}{d\lambda}}\,d\lambda\label{eqIso.4}%
\end{equation}
where $\mu\absolutelycontinuous\lambda$, $\nu\absolutelycontinuous\lambda$.

\item \label{LemIso.2proof(4)}It is known, see \cite{Nel69}, that
$\mathcal{H}$ is a Hilbert space. In particular, it is complete: if a sequence
$\left(  \varphi_{n},\mu_{n}\right)  $ in $\mathcal{H}$ satisfies
\[
\lim_{n,m\rightarrow\infty}\left\Vert \left(  \varphi_{n},\mu_{n}\right)
-\left(  \varphi_{m},\mu_{m}\right)  \right\Vert ^{2}=0,
\]
then there is a pair $\left(  \varphi,\mu\right)  $ with%
\begin{equation}
\int_{X}\left\vert \varphi\right\vert ^{2}\,\frac{d\mu}{d\lambda}%
\,d\lambda=\int_{X}\left\vert \varphi\right\vert ^{2}\,d\mu<\infty,
\label{eqIso.5}%
\end{equation}
where%
\begin{equation}
\lambda:=\sum_{n=1}^{\infty}2^{-n}\mu_{n}\left(  X\right)  ^{-1}\mu_{n},
\label{eqIso.6}%
\end{equation}
and $\left\Vert \left(  \varphi,\mu\right)  -\left(  \varphi_{n},\mu
_{n}\right)  \right\Vert ^{2}\underset{n\rightarrow\infty}{\longrightarrow}0$.
\end{enumerate}

Assuming that the expression in (\ref{eqIso.2}) defines an operator $S$ in
$\mathcal{H}$, it follows from (\ref{eqIso.3}) that $S$ is linear. To see
this, let $\left(  \varphi,\mu\right)  $, $\left(  \psi,\nu\right)  $, and
$\lambda$ be as stated in the conditions below (\ref{eqIso.3}). Then $\mu
\circ\tau^{-1}\absolutelycontinuous\lambda\circ\tau^{-1}$, and $\nu\circ
\tau^{-1}\absolutelycontinuous\lambda\circ\tau^{-1}$, and a calculation shows
that the following formula holds for the transformation of the Radon-Nikodym
derivatives: setting%
\begin{equation}
\frac{d\mu\circ\tau^{-1}}{d\lambda\circ\tau^{-1}}=k_{\mu}, \label{eqIso.7}%
\end{equation}
we have%
\begin{equation}
k_{\mu}\circ\tau=\frac{d\mu}{d\lambda}\qquad(\mathrm{a.e.}\,\lambda).
\label{eqIso.8}%
\end{equation}
Similarly $\displaystyle k_{\nu}:=\frac{d\nu\circ\tau^{-1}}{d\lambda\circ
\tau^{-1}}$ satisfies%
\begin{equation}
k_{\nu}\circ\tau=\frac{d\nu}{d\lambda}\qquad(\mathrm{a.e.}\,\lambda).
\label{eqIso.8half}%
\end{equation}

To show that $S$ is linear, we must calculate the sum%
\begin{equation}
\left(  \varphi\circ\sigma,\mu\circ\tau^{-1}\right)  +\left(  \psi\circ
\sigma,\nu\circ\tau^{-1}\right)  , \label{eqIso.9}%
\end{equation}
or, in expanded notation, we must verify that%
\begin{equation}
\left(  \varphi\circ\sigma\;\sqrt{k_{\mu}}+\psi\circ\sigma\;\sqrt{k_{\nu}%
},\;\lambda\circ\tau^{-1}\right)  \sim\left(  \left(  \varphi\,\sqrt
{\frac{d\mu}{d\lambda}}+\psi\,\sqrt{\frac{d\nu}{d\lambda}}\,\right)
\circ\sigma,\;\lambda\circ\tau^{-1}\right)  _{\mathstrut}. \label{eqIso.10}%
\end{equation}
We get this class identity by an application of (\ref{eqIso.8}) as follows:%
\[
k_{\mu}\left(  x\right)  =k_{\mu}\left(  \tau\left(  \sigma\left(  x\right)
\right)  \right)  =\left.  \left(  \sqrt{\frac{d\mu}{d\lambda}}\circ
\sigma\right)  \right\vert _{\tau\left(  X\right)  }\left(  x\right)
\qquad(\mathrm{a.e.}\,\lambda\circ\tau^{-1}).
\]
Similarly, for the other measure, we get%
\begin{equation}
k_{\nu}=\left.  \left(  \sqrt{\frac{d\nu}{d\lambda}}\circ\sigma\right)
\right\vert _{\tau\left(  X\right)  }\qquad(\mathrm{a.e.}\,\lambda\circ
\tau^{-1}). \label{eqIso.11}%
\end{equation}

Assuming again that $S$ in (\ref{eqIso.2}) is well defined, we now show that
it is isometric, i.e., that $\left\Vert S\left(  \varphi,\mu\right)
\right\Vert ^{2}=\left\Vert \left(  \varphi,\mu\right)  \right\Vert ^{2}$,
referring to the norm of $\mathcal{H}$. In view of (\ref{eqIso.3}) and
(\ref{eqIso.10}), it is enough to show that%
\begin{equation}
\int_{X}\left\vert \varphi\circ\sigma\right\vert ^{2}\,k_{\mu}\;d\lambda
\circ\tau^{-1}=\int_{X}\left\vert \varphi\right\vert ^{2}\,\frac{d\mu
}{d\lambda}\,d\lambda. \label{eqIso.12}%
\end{equation}
But, using (\ref{eqIso.8}), we get
\[
\int_{X}\left\vert \varphi\circ\sigma\right\vert ^{2}\,k_{\mu}\;d\lambda
\circ\tau^{-1}=\int_{X}\left\vert \varphi\circ\sigma\circ\tau\right\vert
^{2}\;k_{\mu}\circ\tau\;d\lambda\underset{\text{(\ref{eqIso.8})}}{=}\int
_{X}\left\vert \varphi\right\vert ^{2}\,\frac{d\mu}{d\lambda}\,d\lambda,
\]
which is the desired formula (\ref{eqIso.12}).

It remains to prove that $S$ is well defined, i.e., that the following
implication holds:%
\begin{equation}
\left(  \varphi,\mu\right)  \sim\left(  \psi,\nu\right)  \Longrightarrow
\left(  \varphi\circ\sigma,\mu\circ\tau^{-1}\right)  \sim\left(  \psi
\circ\sigma,\nu\circ\tau^{-1}\right)  . \label{eqIso.13}%
\end{equation}
To do this, we go through a sequence of implications which again uses the
fundamental transformation rules (\ref{eqIso.8}) and (\ref{eqIso.11}).

Pick some $\lambda$ such that $\mu\absolutelycontinuous\lambda$ and
$\nu\absolutelycontinuous\lambda$. We establish the following implication:
\begin{equation}
\varphi\,\sqrt{\frac{d\mu}{d\lambda}}=\psi\,\sqrt{\frac{d\nu}{d\lambda}%
}\;\left(  \mathrm{a.e.}\,\lambda\right)  \Longrightarrow\left(  \varphi
\circ\sigma\right)  \sqrt{k_{\mu}}=\left(  \psi\circ\sigma\right)
\sqrt{k_{\nu}}\;(\mathrm{a.e.}\,\lambda\circ\tau^{-1}), \label{eqIso.14}%
\end{equation}
where $\displaystyle k_{\mu}=\frac{d\mu\circ\tau^{-1}}{d\lambda\circ\tau^{-1}%
}$ and $\displaystyle k_{\nu}=\frac{d\nu\circ\tau^{-1}}{d\lambda\circ\tau
^{-1}}$. The desired conclusion (\ref{eqIso.13}) follows from this.

We now turn to the proof of the implication (\ref{eqIso.14}). We pick a third
measure $\lambda$ as described, and assume the identity%
\[
\varphi\,\sqrt{\frac{d\mu}{d\lambda}}=\psi\,\sqrt{\frac{d\nu}{d\lambda}}%
\qquad\mathrm{a.e.}\,\lambda.
\]
Let $f$ be a bounded Borel function on $X$. In the following calculations, all
integrals are over the full space $X$, but the measures change as we make
transformations, and we use the definition of the Radon-Nikodym formula. First
note that%
\begin{align*}
\int f\;k_{\mu}\,\left(  \frac{d\nu}{d\lambda}\circ\sigma\right)
\;d\lambda\circ\tau^{-1}  &  =\int f\,\left(  \frac{d\nu}{d\lambda}\circ
\sigma\right)  \;d\mu\circ\tau^{-1}\\
&  =\int f\circ\tau\;\,\frac{d\nu}{d\lambda}\;\,d\mu=\int f\circ\tau
\;\,\frac{d\nu}{d\lambda}\;\,\frac{d\mu}{d\lambda}\;\,d\lambda.
\end{align*}
But by symmetry, we also have%
\[
\int f\;k_{\nu}\,\left(  \frac{d\mu}{d\lambda}\circ\sigma\right)
\;d\lambda\circ\tau^{-1}=\int f\circ\tau\;\frac{d\nu}{d\lambda}\;\frac{d\mu
}{d\lambda}\;d\lambda.
\]
Putting the last two formulas together, we arrive at the following identity:%
\[
\int_{X}f\;k_{\mu}\;\,\frac{d\nu}{d\lambda}\circ\sigma\;\,d\lambda\circ
\tau^{-1}=\int_{X}f\;k_{\nu}\;\,\frac{d\mu}{d\lambda}\circ\sigma
\;\,d\lambda\circ\tau^{-1}.
\]
Since the function $f$ is arbitrary, we get
\[
k_{\mu}\left(  \frac{d\nu}{d\lambda}\circ\sigma\right)  =k_{\nu}\left(
\frac{d\mu}{d\lambda}\circ\sigma\right)  \qquad\mathrm{a.e.}\,\lambda\circ
\tau^{-1}%
\]
and, of course,%
\[
\sqrt{k_{\mu}}\;\,\sqrt{\frac{d\nu}{d\lambda}}\circ\sigma=\sqrt{k_{\nu}%
}\;\,\sqrt{\frac{d\mu}{d\lambda}}\circ\sigma\qquad\mathrm{a.e.}\,\lambda
\circ\tau^{-1}.
\]
Using now the identity%
\[
\varphi\,\sqrt{\frac{d\mu}{d\lambda}}=\psi\,\sqrt{\frac{d\nu}{d\lambda}}%
\qquad\mathrm{a.e.}\,\lambda,
\]
we arrive at the formula%
\[
\varphi\circ\sigma\;\,\sqrt{k_{\mu}}\;\,\sqrt{\frac{d\mu}{d\lambda}}%
\circ\sigma\;\,\sqrt{\frac{d\nu}{d\lambda}}\circ\sigma=\psi\circ
\sigma\;\,\sqrt{k_{\nu}}\;\,\sqrt{\frac{d\mu}{d\lambda}}\circ\sigma
\;\,\sqrt{\frac{d\nu}{d\lambda}}\circ\sigma,
\]
and by cancellation,%
\[
\varphi\circ\sigma\;\sqrt{k_{\mu}}=\psi\circ\sigma\;\sqrt{k_{\nu}}%
\qquad\mathrm{a.e.}\,\lambda\circ\tau^{-1}.
\]
This completes the proof of the implication (\ref{eqIso.14}), and therefore
also of (\ref{eqIso.13}). This means that if the linear operator $S$ is
defined as in (\ref{eqIso.2}), then the result is independent of which element
is chosen in the equivalence class represented by the pair $\left(
\varphi,\mu\right)  $. Putting together the steps in the proof, we conclude
that $S\colon\mathcal{H}\rightarrow\mathcal{H}$ is an isometry, and that it
has the properties which are stated in the lemma. \qed
\end{pf}

\begin{lem}
\label{LemIso.3}Let $X$ be a compact Hausdorff space, and let $\sigma$ be as
in the statement of Lemma \textup{\ref{LemIso.2}}, i.e., $\sigma\colon
X\rightarrow X$ is onto and continuous. Suppose $\sigma$ has two distinct
branches of the inverse, i.e., $\tau_{i}\colon X\rightarrow X$, $i=1,2$,
continuous, and satisfying $\sigma\circ\tau_{i}=\operatorname*{id}%
\nolimits_{X}$, $i=1,2$. Let $S_{i}\colon\mathcal{H}\rightarrow\mathcal{H}$ be
the corresponding isometries, i.e.,%
\begin{equation}
S_{i}\left(  \left(  \varphi,\mu\right)  \right)  :=\left(  \varphi\circ
\sigma,\mu\circ\tau_{i}^{-1}\right)  , \label{eqIso.15}%
\end{equation}
or%
\[
S_{i}\colon\varphi\sqrt{d\mu}\longmapsto\varphi\circ\sigma\;\sqrt{d\mu
\circ\tau_{i}^{-1}}.
\]
Then the two isometries have orthogonal ranges, i.e.,%
\begin{equation}
\left\langle \,S_{1}\left(  \left(  \varphi,\mu\right)  \right)  \mid
S_{2}\left(  \left(  \psi,\nu\right)  \right)  \right\rangle =0
\label{eqIso.16}%
\end{equation}
for all pairs of vectors in $\mathcal{H}$, i.e., all $\left(  \varphi
,\mu\right)  \in\mathcal{H}$ and $\left(  \psi,\nu\right)  \in\mathcal{H}$.
\end{lem}

\begin{pf}
Note that in the statement (\ref{eqIso.16}) of the conclusion, we use
$\left\langle \,\cdot\mid\cdot\,\right\rangle $ to denote the inner product of
the Hilbert space $\mathcal{H}$, as it was defined in (\ref{eqIso.4}).

With the two measures $\mu$ and $\nu$ given, then the expression in
(\ref{eqIso.16}) involves the transformed measures $\mu\circ\tau_{1}^{-1}$ and
$\nu\circ\tau_{2}^{-1}$. Now pick some measure $\lambda$ such that $\mu
\circ\tau_{1}^{-1}\absolutelycontinuous\lambda$ and $\nu\circ\tau_{2}%
^{-1}\absolutelycontinuous\lambda$. Then the expression in (\ref{eqIso.16}) is%
\begin{equation}
\int_{X}\overline{\varphi\circ\sigma}\;\;\psi\circ\sigma\;\,\sqrt{\frac
{d\mu\circ\tau_{1}^{-1}}{d\lambda}}\;\,\sqrt{\frac{d\nu\circ\tau_{2}^{-1}%
}{d\lambda}}\;\,d\lambda. \label{eqIso.17}%
\end{equation}
But $\displaystyle\frac{d\mu\circ\tau_{1}^{-1}}{d\lambda}$ is supported on
$\tau_{1}\left(  X\right)  $, while $\displaystyle\frac{d\nu\circ\tau_{2}%
^{-1}}{d\lambda}$ is supported on $\tau_{2}\left(  X\right)  $. Since
$\tau_{1}\left(  X\right)  \cap\tau_{2}\left(  X\right)  ^{\mathstrut
}=\varnothing$ by the choice of distinct branches for the inverse of $\sigma$,
we conclude that the integral in (\ref{eqIso.17}) vanishes. \qed
\end{pf}

In the next lemma we prove a formula for the adjoint $S^{\ast}$ of the
isometry $S$ which was introduced in Lemma \ref{LemIso.2}. Now $S^{\ast}$
refers to the inner product (\ref{eqIso.4}) which is given at the outset, and
which defines the Hilbert space $\mathcal{H}$.

\begin{lem}
\label{LemIso.4}Let $X$, $\sigma$, and $\tau$ be given as in the statement of
Lemma \textup{\ref{LemIso.2}}, i.e., we assume that $\sigma$ is onto, that
both $\sigma$ and $\tau$ are continuous, and that
\begin{equation}
\sigma\circ\tau=\operatorname*{id}\nolimits_{X}.\label{eqIso.18}%
\end{equation}
Let $S$ be the isometry defined in \textup{(\ref{eqIso.2})}, and let $S^{\ast
}$ be the adjoint co-isometry. Then
\begin{equation}
S^{\ast}\left(  \left(  \varphi,\mu\right)  \right)  =\left(  \varphi\circ
\tau,\mu\circ\sigma^{-1}\right)  \label{eqIso.19}%
\end{equation}
for all $\left(  \varphi,\mu\right)  \in S\mkern2mu\mathcal{H}$.
\end{lem}

\begin{pf}
Recall that operators in $\mathcal{H}$ are defined on equivalence classes:
just as in the proof of Lemma \ref{LemIso.2}, we must check the implication%
\begin{equation}
\left(  \varphi,\mu\right)  \sim\left(  \psi,\nu\right)  \Longrightarrow
\left(  \varphi\circ\tau,\mu\circ\sigma^{-1}\right)  \sim\left(  \psi\circ
\tau,\nu\circ\sigma^{-1}\right)  . \label{eqIso.20}%
\end{equation}
While the verification of (\ref{eqIso.20}) involves the transformation rules
for Radon-Nikodym derivatives, the steps of the proof are quite analogous to
the arguments from the proof of Lemma \ref{LemIso.2}, and they are left to the reader.

Now let $T$ denote the operator on $\mathcal{H}$ which is defined by the
formula (\ref{eqIso.19}). It is clear that $TS=I={}$the identity operator in
$\mathcal{H}$, i.e., that $TS\left(  \varphi,\mu\right)  =\left(  \varphi
,\mu\right)  $ for all $\left(  \varphi,\mu\right)  \in\mathcal{H}$. Indeed,%
\[
TS\left(  \varphi,\mu\right)  =T\left(  \varphi\circ\sigma,\;\mu\circ\tau
^{-1}\right)  =\left(  \varphi\circ\sigma\circ\tau,\;\mu\circ\tau^{-1}%
\circ\sigma^{-1}\right)  =\left(  \varphi,\mu\right)  ,
\]
where the identity (\ref{eqIso.18}) was used in the last step of the argument.

The assertion of the lemma is that $T|_{S\mathcal{H}}=S^{\ast}$. Since $TS=I$,
and $S$ is isometric, we need only set $T$ equal to zero on $\left(
S\mathcal{H}\right)  ^{\perp}$, where
\begin{equation}
\left(  S\mathcal{H}\right)  ^{\perp}=\left\{  \,x\in\mathcal{H}%
\mid\left\langle \,Sy\mid x\,\right\rangle =0,\;y\in\mathcal{H}\,\right\}  .
\label{eqIso.21}%
\end{equation}
Let $x=\operatorname*{class}\left(  \psi,\nu\right)  \in\left(  S\mathcal{H}%
\right)  ^{\perp}$, and let $y=\operatorname*{class}\left(  \varphi
,\mu\right)  $. The argument from the proof of Lemma \ref{LemIso.2} shows that
there is a positive Borel measure $\lambda$ such that $\mu\circ\tau
^{-1}\absolutelycontinuous\lambda\circ\tau^{-1}$ and $\nu
\absolutelycontinuous\lambda\circ\tau^{-1}$. Recall that the Radon-Nikodym
derivative $\displaystyle k_{\mu}=\frac{d\mu\circ\tau^{-1}}{d\lambda\circ
\tau^{-1}}$ satisfies
\begin{equation}
k_{\mu}\circ\tau=\frac{d\mu}{d\lambda}. \label{eqIso.22}%
\end{equation}
We now calculate the inner-product term from (\ref{eqIso.21}):%
\begin{align*}
\left\langle \,Sy\mid x\,\right\rangle  &  =\int_{X}\overline{\varphi
\circ\sigma}\;\;\psi\;\,\sqrt{k_{\mu}}\;\,\sqrt{\frac{d\nu}{d\lambda\circ
\tau^{-1}}}\;\,d\lambda\circ\tau^{-1}\\[3\jot]
&  =\int_{X}\bar{\varphi}\;\,\psi\circ\tau\;\,\sqrt{k_{\mu}}\circ\tau
\;\,\sqrt{\frac{d\nu}{d\lambda\circ\tau^{-1}}}\circ\tau\;\,d\lambda\\[3\jot]
&  =\int_{X}\bar{\varphi}\;\,\psi\circ\tau\;\,\sqrt{\frac{d\mu}{d\lambda}%
}\;\,\sqrt{\frac{d\nu}{d\lambda\circ\tau^{-1}}}\circ\tau\;\,d\lambda,
\end{align*}
where we used (\ref{eqIso.22}) in the last step.

Since this expression${}\equiv0$ for all $\left(  \varphi,\mu\right)
\in\mathcal{H}$, and $\displaystyle d\mu=\frac{d\mu}{d\lambda}\,d\lambda$, we
conclude that $\psi\circ\tau=0$\enspace$\mathrm{a.e.}\,\lambda$. \qed
\end{pf}

\section{\label{Rep}Representations in $L^{2}\left(  \mu\right)  $}

Let $X$ be a compact Hausdorff space, and let $\mu$ be a positive Borel
measure on $X$. For simplicity, we assume that $\mu$ is normalized, i.e., that
$\mu\left(  X\right)  =1$. Ideally, we look for some measure $\mu$ such that
the Hilbert space $L^{2}\left(  \mu\right)  =L^{2}\left(  X,\mu\right)  $
suffices for the representation theory under discussion. For representations
of the Cuntz algebras $\mathcal{O}_{N}$, it may be possible to stay within the
Hilbert space $L^{2}\left(  \mu\right)  $, while for some other
representations, the \textquotedblleft larger\textquotedblright\ Hilbert space
$\mathcal{H}$ of Section \ref{Iso} is forced on us.

The general setting in this section will be the same as in Section \ref{Iso}:
the transformations $\sigma,\tau\colon X\rightarrow X$ are assumed continuous,
and $\tau$ is a branch of the inverse of $\sigma$, i.e., we assume that%
\begin{equation}
\sigma\circ\tau=\operatorname*{id}\nolimits_{X}. \label{eqRep.1}%
\end{equation}
It follows that $\sigma$ is onto, and that $\tau$ is one-to-one. We will show
in this section that if%
\begin{equation}
\mu\circ\tau^{-1}\absolutelycontinuous\mu, \label{eqRep.pound}%
\end{equation}
then the isometry $S$ from Lemma \ref{LemIso.2} may be realized in
$L^{2}\left(  \mu\right)  $.

\begin{thm}
\label{ThmRep.1}Let $\sigma,\tau\colon X\rightarrow X$ be continuous, and
suppose that $\sigma\circ\tau=\operatorname*{id}\nolimits_{X}$ holds. Let
$\mu$ be a positive Borel measure on $X$ such that $\mu\left(  X\right)  =1$
and%
\begin{equation}
\mu\circ\tau^{-1}\absolutelycontinuous\mu. \label{eqRep.2}%
\end{equation}
Let $S\colon\mathcal{H}\rightarrow\mathcal{H}$ be the isometry defined in
Lemma \textup{\ref{LemIso.2}}, i.e.,%
\begin{equation}
S\left(  \left(  \varphi,\mu\right)  \right)  :=\left(  \varphi\circ\sigma
,\mu\circ\tau^{-1}\right)  . \label{eqRep.3}%
\end{equation}
Setting
\begin{equation}
S_{\mu}\varphi:=\varphi\circ\sigma\;\sqrt{\frac{d\mu\circ\tau^{-1}}{d\mu}}
\label{eqRep.4}%
\end{equation}
and%
\begin{equation}
W_{\mu}\varphi:=\left(  \varphi,\mu\right)  , \label{eqRep.5}%
\end{equation}
we get two isometries, $S_{\mu}\colon L^{2}\left(  \mu\right)  \rightarrow
L^{2}\left(  \mu\right)  $ and $W_{\mu}\colon L^{2}\left(  \mu\right)
\rightarrow\mathcal{H}$, such that%
\begin{equation}
W_{\mu}S_{\mu}=SW_{\mu}. \label{eqRep.6}%
\end{equation}

\end{thm}

\begin{figure}[ptb]
\begin{center}
$\displaystyle%
\begin{array}
[c]{ccc}%
L^{2}\left(  \mu\right)  & \overset{W_{\mu}}{\longrightarrow} & \mathcal{H}\\
\llap{$\scriptstyle S_{\mu}$}\downarrow &  & \downarrow
\rlap{$\scriptstyle S$}\\
L^{2}\left(  \mu\right)  & \overset{W_{\mu}}{\longrightarrow} & \mathcal{H}%
\end{array}
$
\end{center}
\caption{$W_{\mu}S_{\mu}=SW_{\mu}$.}%
\label{Fig1}%
\end{figure}

\begin{pf}
The assertion (\ref{eqRep.6}) states that $W_{\mu}$ intertwines the two
isometries $S_{\mu}$ and $S$, or, expressed as a diagram, that the
commutativity shown in Fig.\ \ref{Fig1} holds. To prove that $S_{\mu}$ is
isometric, note that
\begin{align*}
\int_{X}\left\vert S_{\mu}\varphi\right\vert \,d\mu &  =\int_{X}\left\vert
\varphi\circ\sigma\right\vert ^{2}\;\frac{d\mu\circ\tau^{-1}}{d\mu}\;d\mu
=\int_{X}\left\vert \varphi\circ\sigma\right\vert ^{2}\;d\mu\circ\tau^{-1}\\
&  =\int_{X}\left\vert \varphi\circ\sigma\circ\tau\right\vert ^{2}\;d\mu
=\int_{X}\left\vert \varphi\right\vert ^{2}\,d\mu.
\end{align*}

It is clear from the definition of the norm in $\mathcal{H}$ that $W_{\mu}$ is
isometric. To verify (\ref{eqRep.6}), we note that%
\[
W_{\mu}S_{\mu}\varphi=\left(  \varphi\circ\sigma\;\,\sqrt{\frac{d\mu\circ
\tau^{-1}}{d\mu}},\;\mu\right)  ,
\]
and that%
\[
SW_{\mu}\varphi=\left(  \varphi\circ\sigma,\mu\circ\tau^{-1}\right)  .
\]
But since $\mu\circ\tau^{-1}\absolutelycontinuous\mu$, it is clear from
(\ref{eqIso.1}) that%
\[
\left(  \varphi\circ\sigma\;\,\sqrt{\frac{d\mu\circ\tau^{-1}}{d\mu}}%
,\;\mu\right)  \sim\left(  \varphi\circ\sigma,\mu\circ\tau^{-1}\right)  .
\]
Since the vectors in $\mathcal{H}$ are equivalence classes, the desired
intertwining identity (\ref{eqRep.6}) holds. \qed
\end{pf}

\begin{cor}
\label{CorRep.2}Let $X$ be a compact Hausdorff space, and let $\sigma
,\tau\colon X\rightarrow X$ satisfy the conditions stated in Theorem
\textup{\ref{ThmRep.1}}. Let $\mu$ be a positive Borel measure on $X$ such
that $\mu\left(  X\right)  =1$ and $\mu\circ\tau^{-1}\absolutelycontinuous\mu
$. Then the Radon-Nikodym de\-riv\-a\-tive%
\begin{equation}
p_{\mu}:=\frac{d\mu\circ\tau^{-1}}{d\mu} \label{eqRep.7}%
\end{equation}
satisfies%
\begin{align}
p_{\mu}  &  \geq1\qquad\mu\text{-}\mathrm{a.e.}\text{ on }\tau\left(
X\right)  ,\label{eqRep.8}\\[3pt]
S_{\mu}^{\ast}\varphi &  =\varphi\circ\tau\;\left(  p_{\mu}\circ\tau\right)
^{-1/2}. \label{eqRep.9}%
\end{align}

\end{cor}

\begin{pf}
Since $S_{\mu}\colon L^{2}\left(  \mu\right)  \rightarrow L^{2}\left(
\mu\right)  $ is isometric by the theorem, $S_{\mu}^{\ast}$ is contractive in
$L^{2}\left(  \mu\right)  $, and $\left\Vert S_{\mu}^{\ast}\right\Vert =1$.
But a substitution of formula (\ref{eqRep.9}) yields $\left\langle
\,\varphi\mid S_{\mu}\psi\,\right\rangle =\left\langle
\,\smash{S_{\mu }^{\ast }}\varphi\mid\psi\,\right\rangle $ for $\varphi
,\psi\in L^{2}\left(  \mu\right)  $. Indeed, we have the following identity:%
\begin{align*}
\int\bar{\varphi}\;\,\psi\circ\sigma\;\,p_{\mu}^{1/2}\;\,d\mu &  =\int
\bar{\varphi}\;\,\psi\circ\sigma\;\,p_{\mu}^{-1/2}\;p_{\mu}\;\,d\mu\\
&  =\int\bar{\varphi}\;\,\psi\circ\sigma\;\,p_{\mu}^{-1/2}\;\,d\mu\circ
\tau^{-1}=\int\overline{\varphi\circ\tau}\;\,\psi\;\left(  p_{\mu}\circ
\tau\right)  ^{-1/2}\;\,d\mu.
\end{align*}
This proves formula (\ref{eqRep.9}) for the co-isometry $S_{\mu}^{\ast}\colon
L^{2}\left(  \mu\right)  \rightarrow L^{2}\left(  \mu\right)  $. \qed
\end{pf}

\begin{cor}
\label{CorRep.3}Let $X$ be a compact Hausdorff space, and let $N\in\mathbb{N}%
$, $N\geq2$, be given. Let $\sigma\colon X\rightarrow X$ be continuous and
onto. Suppose there are $N$ distinct branches of the inverse, i.e., continuous
$\tau_{i}\colon X\rightarrow X$, $i=1,\dots,N$, such that%
\begin{equation}
\sigma\circ\tau_{i}=\operatorname*{id}\nolimits_{X}. \label{eqRep.9bis}%
\end{equation}
Suppose there is a positive Borel measure $\mu$ such that $\mu\left(
X\right)  =1$, and%
\begin{equation}
\mu\circ\tau_{i}^{-1}\absolutelycontinuous\mu,\qquad i=1,\dots,N.
\label{eqRep.10}%
\end{equation}
Then the isometries%
\begin{equation}
S_{i}\varphi:=\varphi\circ\sigma\;\sqrt{\frac{d\mu\circ\tau_{i}^{-1}}{d\mu}}
\label{eqRep.11}%
\end{equation}
satisfy%
\begin{equation}
\sum_{i=1}^{N}S_{i}S_{i}^{\ast}=I_{L^{2}\left(  \mu\right)  } \label{eqRep.12}%
\end{equation}
if and only if%
\begin{equation}
\bigcup_{i=1}^{N}\tau_{i}\left(  X\right)  =X. \label{eqRep.13}%
\end{equation}

\end{cor}

\begin{pf}
We already know from Lemma \ref{LemIso.3} that the isometries $S_{i}\colon
L^{2}\left(  \mu\right)  \rightarrow L^{2}\left(  \mu\right)  $ are mutually
orthogonal, i.e., that%
\begin{equation}
S_{i}^{\ast}S_{j}=\delta_{i,j}I_{L^{2}\left(  \mu\right)  }. \label{eqRep.14}%
\end{equation}
It follows that the terms in the sum (\ref{eqRep.12}) are commuting
projections. Hence%
\begin{equation}
\sum_{i=1}^{N}S_{i}S_{i}^{\ast}\leq I_{L^{2}\left(  \mu\right)  }.
\label{eqRep.15}%
\end{equation}
Moreover, we conclude that (\ref{eqRep.12}) holds if and only if%
\begin{equation}
\sum_{i=1}^{N}\left\Vert S_{i}^{\ast}\varphi\right\Vert ^{2}=\left\Vert
\varphi\right\Vert ^{2},\qquad\varphi\in L^{2}\left(  \mu\right)  .
\label{eqRep.16}%
\end{equation}
Setting
$\smash[b]{\displaystyle p_{i}:=\frac{d\mu \circ \tau _{i}^{-1}}{d\mu }}$, we
get%
\begin{equation}
S_{i}^{\ast}\varphi=\varphi\circ\tau_{i}\;\left(  p_{i}\circ\tau_{i}\right)
^{-1/2}; \label{eqRep.17}%
\end{equation}
see (\ref{eqRep.9}) of the previous corollary. We get%
\[
\left\Vert S_{i}^{\ast}\varphi\right\Vert ^{2}=\int_{X}\left\vert \varphi
\circ\tau_{i}\right\vert ^{2}\;\left(  p_{i}\circ\tau_{i}\right)  ^{-1}%
\,d\mu=\int_{\tau_{i}\left(  X\right)  }\left\vert \varphi\right\vert
^{2}\;p_{i}^{-1}\;d\mu\circ\tau_{i}^{-1}=\int_{\tau_{i}\left(  X\right)
}\left\vert \varphi\right\vert ^{2}\;d\mu.
\]
Recall that the branches $\tau_{i}$ of the inverse are distinct. So in view of
(\ref{eqRep.9}), the sets $\tau_{i}\left(  X\right)  $ are non-overlapping.
The equivalence (\ref{eqRep.12})${}\Leftrightarrow{}$(\ref{eqRep.13}) now
follows directly from the previous calculation. \qed
\end{pf}

\begin{figure}[ptb]
\begin{center}
\setlength{\unitlength}{0.3bp} \begin{picture}(390,450)(-15,-75)
\put(0,1){\includegraphics[bb=0 0 360 360,width=108bp,height=108bp]{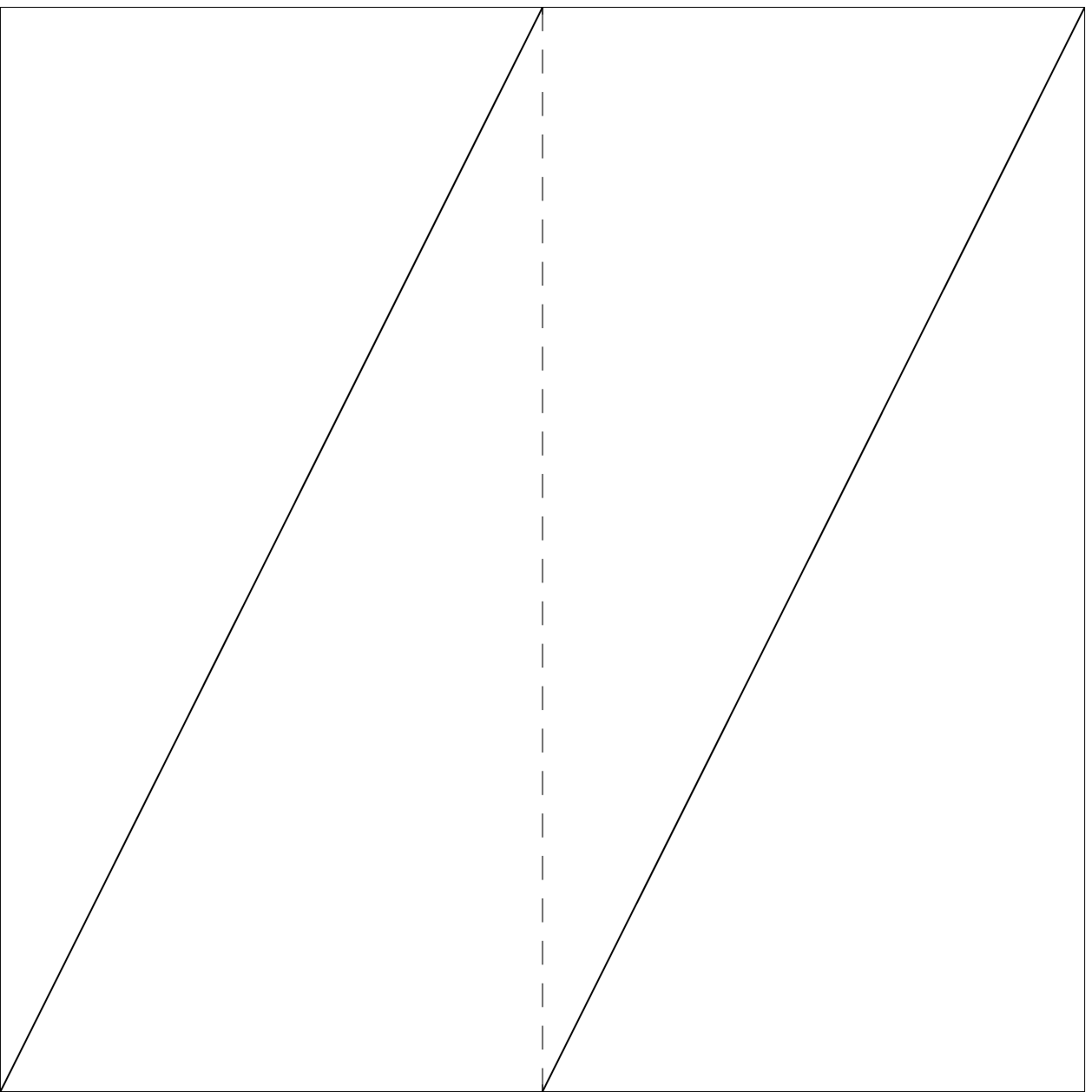}}
\put(-2,-2){\makebox(0,0)[tr]{$0$}}
\put(180,-2){\makebox(0,0)[t]{$\frac{1}{2}$}}
\put(360,-2){\makebox(0,0)[t]{$1$}}
\put(-2,360){\makebox(0,0)[r]{$1$}}
\put(180,-62){\makebox(0,0)[t]{$\sigma$}}
\end{picture}\kern60\unitlength\begin{picture}(390,450)(-15,-75)
\put(0,1){\includegraphics[bb=0 0 360 360,width=108bp,height=108bp]{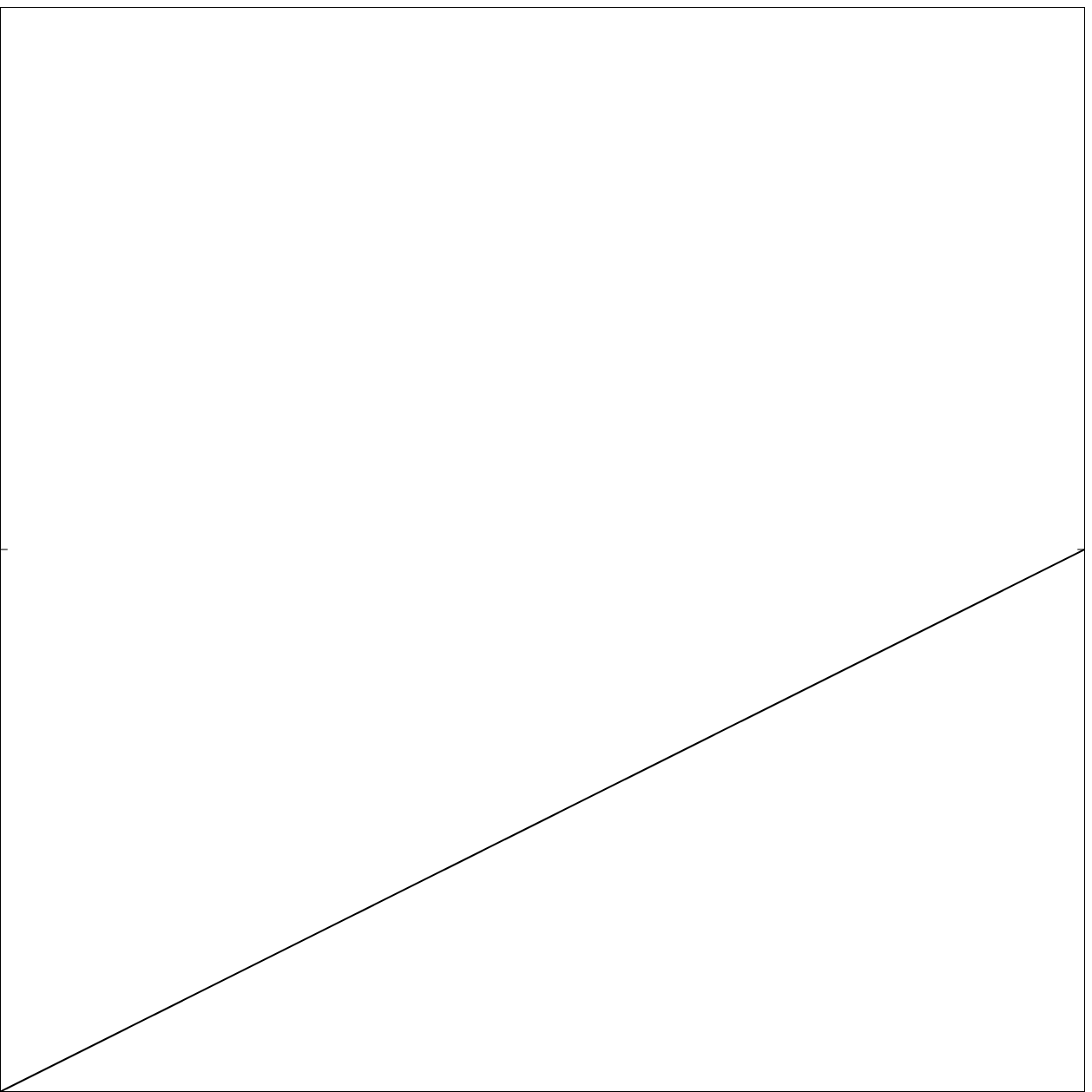}}
\put(-2,-2){\makebox(0,0)[tr]{$0$}}
\put(360,-2){\makebox(0,0)[t]{$1$}}
\put(-2,360){\makebox(0,0)[r]{$1$}}
\put(-2,180){\makebox(0,0)[r]{$\frac{1}{2}$}}
\put(180,-62){\makebox(0,0)[t]{$\tau_{1}$}}
\end{picture}\kern60\unitlength\begin{picture}(390,450)(-15,-75)
\put(0,1){\includegraphics[bb=0 0 360 360,width=108bp,height=108bp]{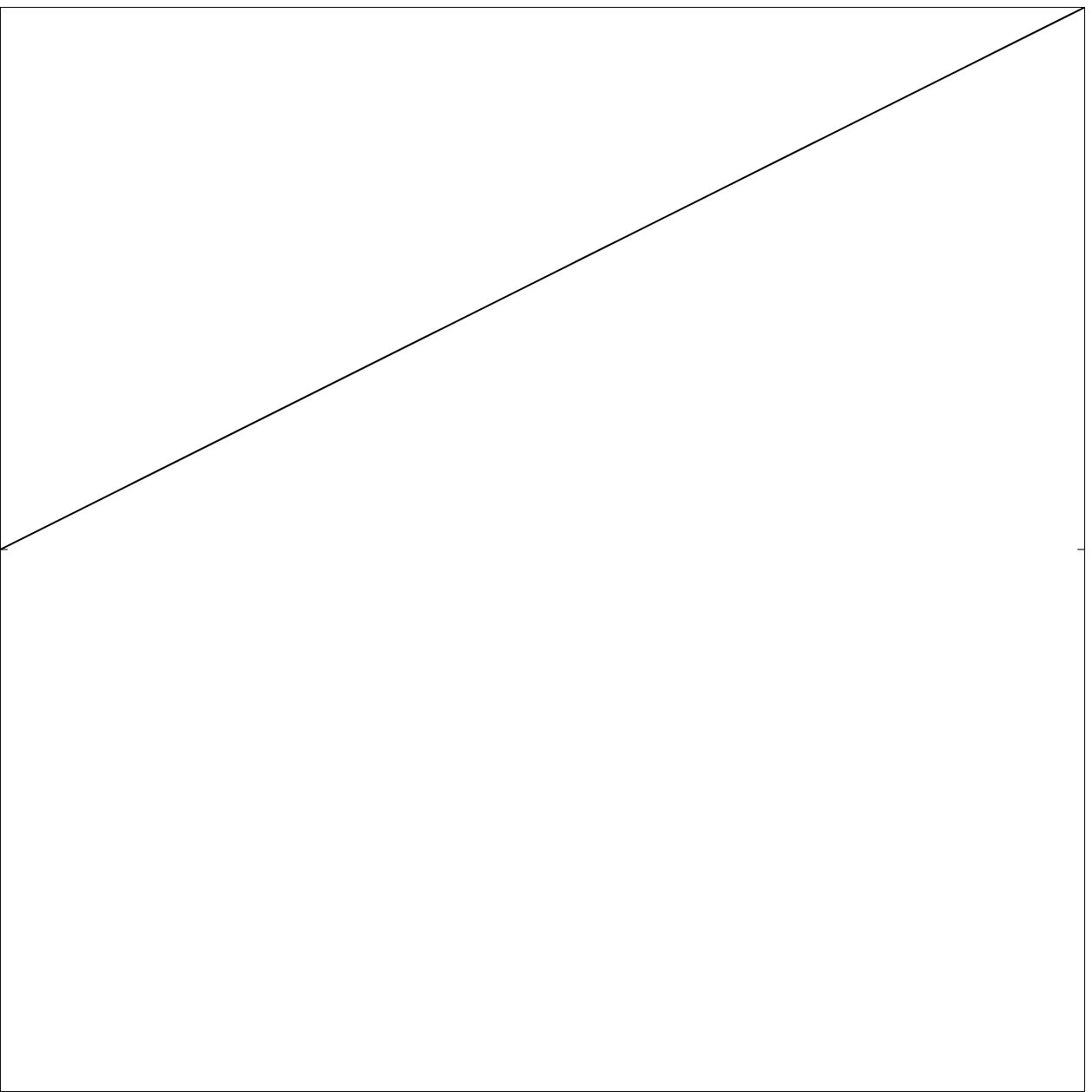}}
\put(-2,-2){\makebox(0,0)[tr]{$0$}}
\put(360,-2){\makebox(0,0)[t]{$1$}}
\put(-2,360){\makebox(0,0)[r]{$1$}}
\put(-2,180){\makebox(0,0)[r]{$\frac{1}{2}$}}
\put(180,-62){\makebox(0,0)[t]{$\tau_{2}$}}
\end{picture}
\end{center}
\caption{Subdivisions of the unit interval.}%
\label{Fig2}%
\end{figure}

\begin{exmp}
\label{ExaRep.4}\upshape Let $X$ be the unit interval $\left[  0,1\right)  $
and let $\mu$ be the restricted Lebesgue measure. Let $\sigma\left(  x\right)
=2x\operatorname{mod}1$, and set $\tau_{1}\left(  x\right)  =\frac{x}{2}$,
$\tau_{2}\left(  x\right)  =\frac{x+1}{2}$. The graphs of the three maps are
illustrated in Figure \textup{\ref{Fig2}}. Then the two isometries%
\begin{equation}
S_{1}\varphi\left(  x\right)  =\varphi\left(  \sigma\left(  x\right)  \right)
\sqrt{2}\chi_{\left[  0,1/2\right)  }\left(  x\right)  ,\;S_{2}\varphi\left(
x\right)  =\varphi\left(  \sigma\left(  x\right)  \right)  \sqrt{2}%
\chi_{\left[  1/2,1\right)  }\left(  x\right)  , \label{eqRep.18}%
\end{equation}
and their adjoints%
\begin{equation}
S_{1}^{\ast}\varphi\left(  x\right)  =\frac{1}{\sqrt{2}}\varphi\left(
\frac{x}{2}\right)  ,\quad S_{2}^{\ast}\varphi\left(  x\right)  =\frac
{1}{\sqrt{2}}\varphi\left(  \frac{x+1}{2}\right)  , \label{eqRep.19}%
\end{equation}
satisfy the relations (\ref{eqRep.12}) and (\ref{eqRep.14}), i.e., they define
a representation of the Cuntz algebra $\mathcal{O}_{2}$ on the Hilbert space
$L^{2}\left(  0,1\right)  $. Note that the conditions (\ref{eqRep.10}) are
satisfied since%
\begin{equation}
\mu\circ\tau_{1}^{-1}=2\mu|_{\left[  0,1/2\right)  }\text{\qquad and\qquad}%
\mu\circ\tau_{2}^{-1}=2\mu|_{\left[  1/2,1\right)  }. \label{eqRep.20}%
\end{equation}

\end{exmp}

\begin{figure}[ptb]
\begin{center}
\setlength{\unitlength}{0.3bp} \begin{picture}(390,555)(-15,-180)
\put(0,1){\includegraphics[bb=0 0 360 360,width=108bp,height=108bp]{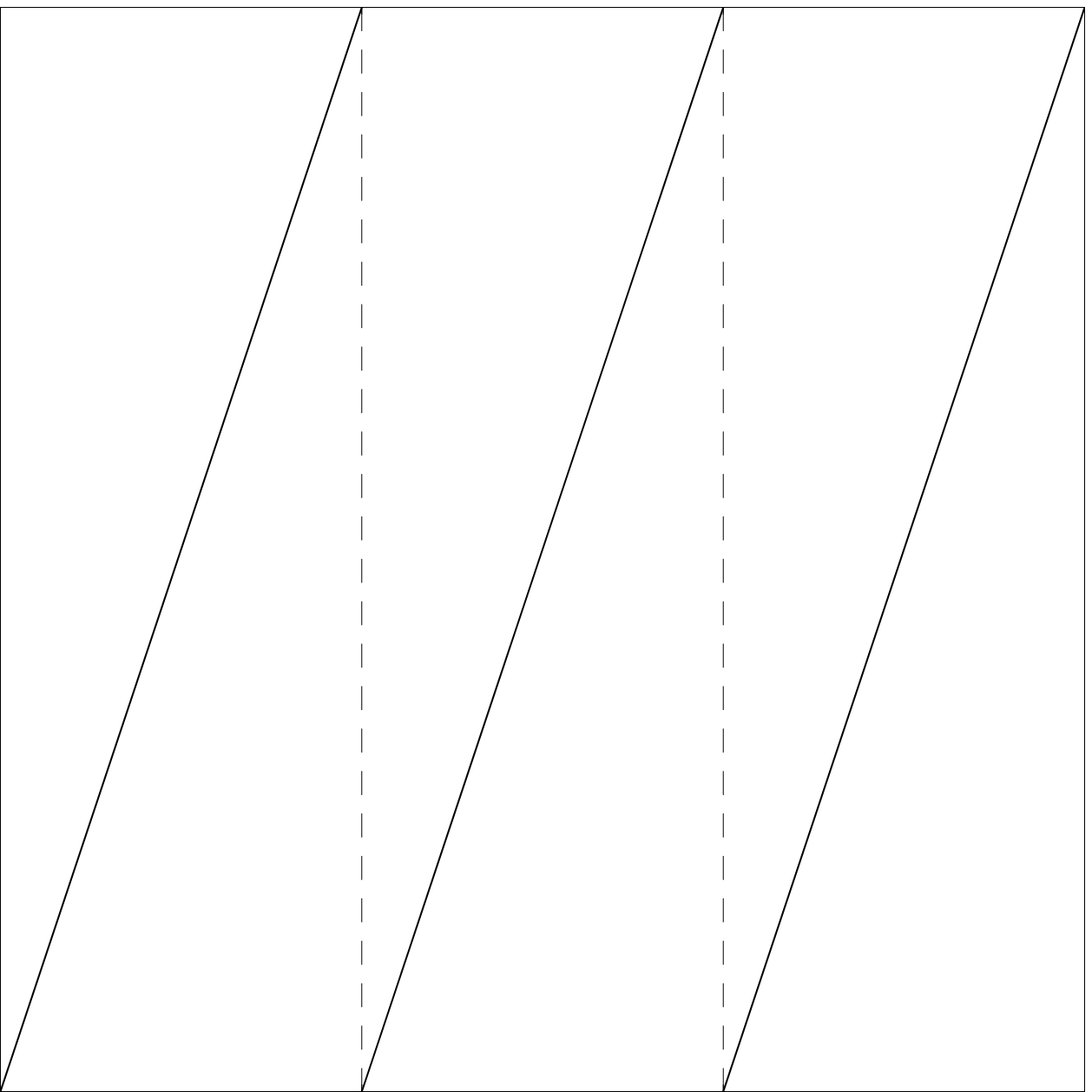}}
\put(-2,-2){\makebox(0,0)[tr]{$0$}}
\put(120,-2){\makebox(0,0)[t]{$\frac{1}{3}$}}
\put(240,-2){\makebox(0,0)[t]{$\frac{2}{3}$}}
\put(360,-2){\makebox(0,0)[t]{$1$}}
\put(-2,360){\makebox(0,0)[r]{$1$}}
\put(0,-69){\includegraphics[bb=4 0 356 13,width=108bp,height=3bp]{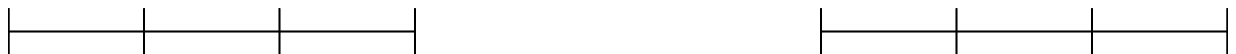}}
\put(0,-88){\includegraphics[bb=4 0 356 13,width=36bp,height=3bp]{cantline.eps}}
\put(240,-88){\includegraphics[bb=4 0 356 13,width=36bp,height=3bp]{cantline.eps}}
\put(0,-107){\includegraphics[bb=4 0 356 13,width=12bp,height=3bp]{cantline.eps}}
\put(80,-107){\includegraphics[bb=4 0 356 13,width=12bp,height=3bp]{cantline.eps}}
\put(240,-107){\includegraphics[bb=4 0 356 13,width=12bp,height=3bp]{cantline.eps}}
\put(320,-107){\includegraphics[bb=4 0 356 13,width=12bp,height=3bp]{cantline.eps}}
\put(6.7,-102){\makebox(0,0)[t]{$\vdots$}}
\put(33.3,-102){\makebox(0,0)[t]{$\vdots$}}
\put(86.7,-102){\makebox(0,0)[t]{$\vdots$}}
\put(113.3,-102){\makebox(0,0)[t]{$\vdots$}}
\put(246.7,-102){\makebox(0,0)[t]{$\vdots$}}
\put(273.3,-102){\makebox(0,0)[t]{$\vdots$}}
\put(326.7,-102){\makebox(0,0)[t]{$\vdots$}}
\put(353.3,-102){\makebox(0,0)[t]{$\vdots$}}
\put(180,-153){\makebox(0,0)[t]{$\underbrace{\rule{108bp}{0bp}}_{\textstyle\sigma}$}}
\end{picture}\kern60\unitlength\begin{picture}(390,555)(-15,-180)
\put(0,1){\includegraphics[bb=0 0 360 360,width=108bp,height=108bp]{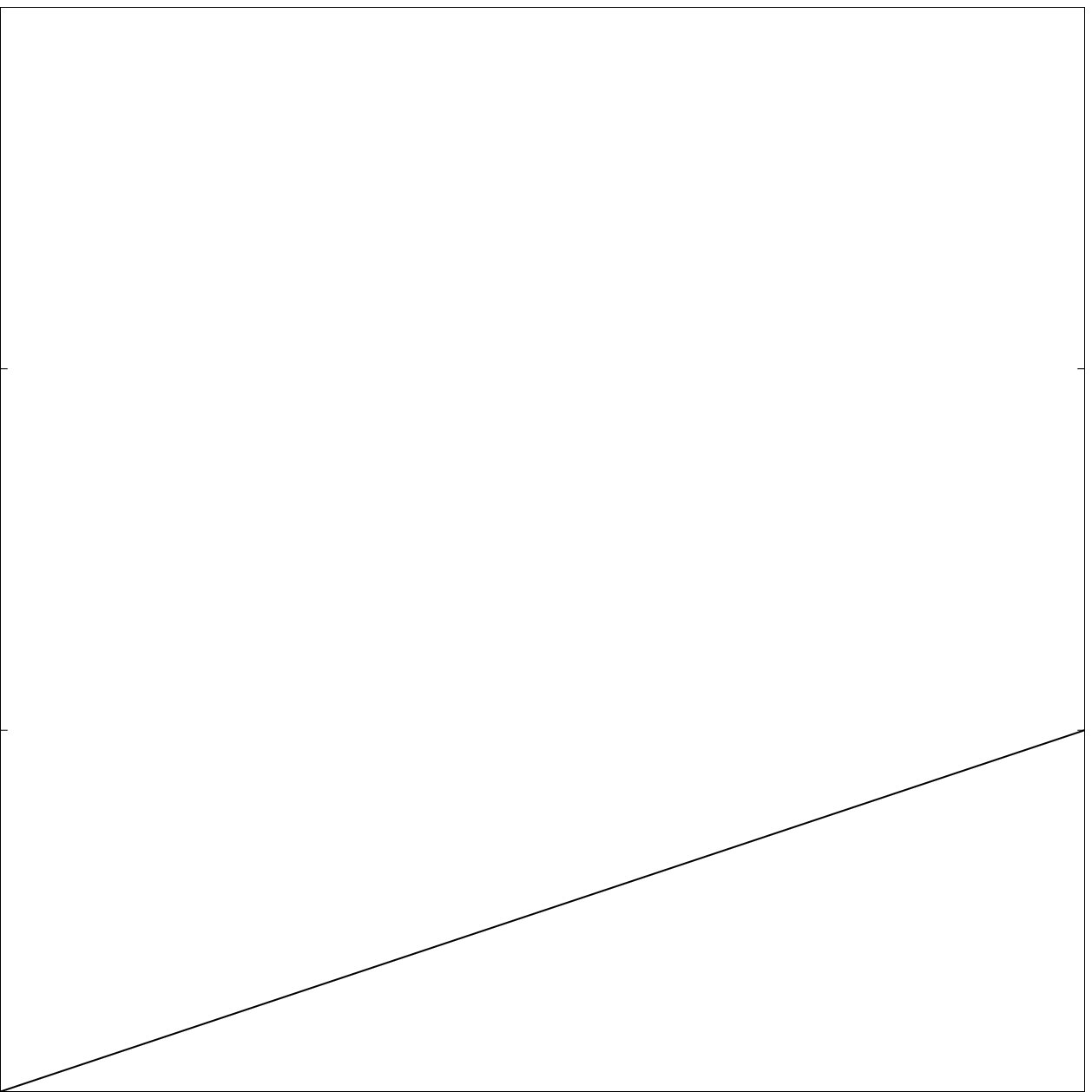}}
\put(-2,-2){\makebox(0,0)[tr]{$0$}}
\put(360,-2){\makebox(0,0)[t]{$1$}}
\put(-2,360){\makebox(0,0)[r]{$1$}}
\put(-2,120){\makebox(0,0)[r]{$\frac{1}{3}$}}
\put(-2,240){\makebox(0,0)[r]{$\frac{2}{3}$}}
\put(180,-62){\makebox(0,0)[t]{$\tau_{1}$}}
\end{picture}\kern60\unitlength\begin{picture}(390,555)(-15,-180)
\put(0,1){\includegraphics[bb=0 0 360 360,width=108bp,height=108bp]{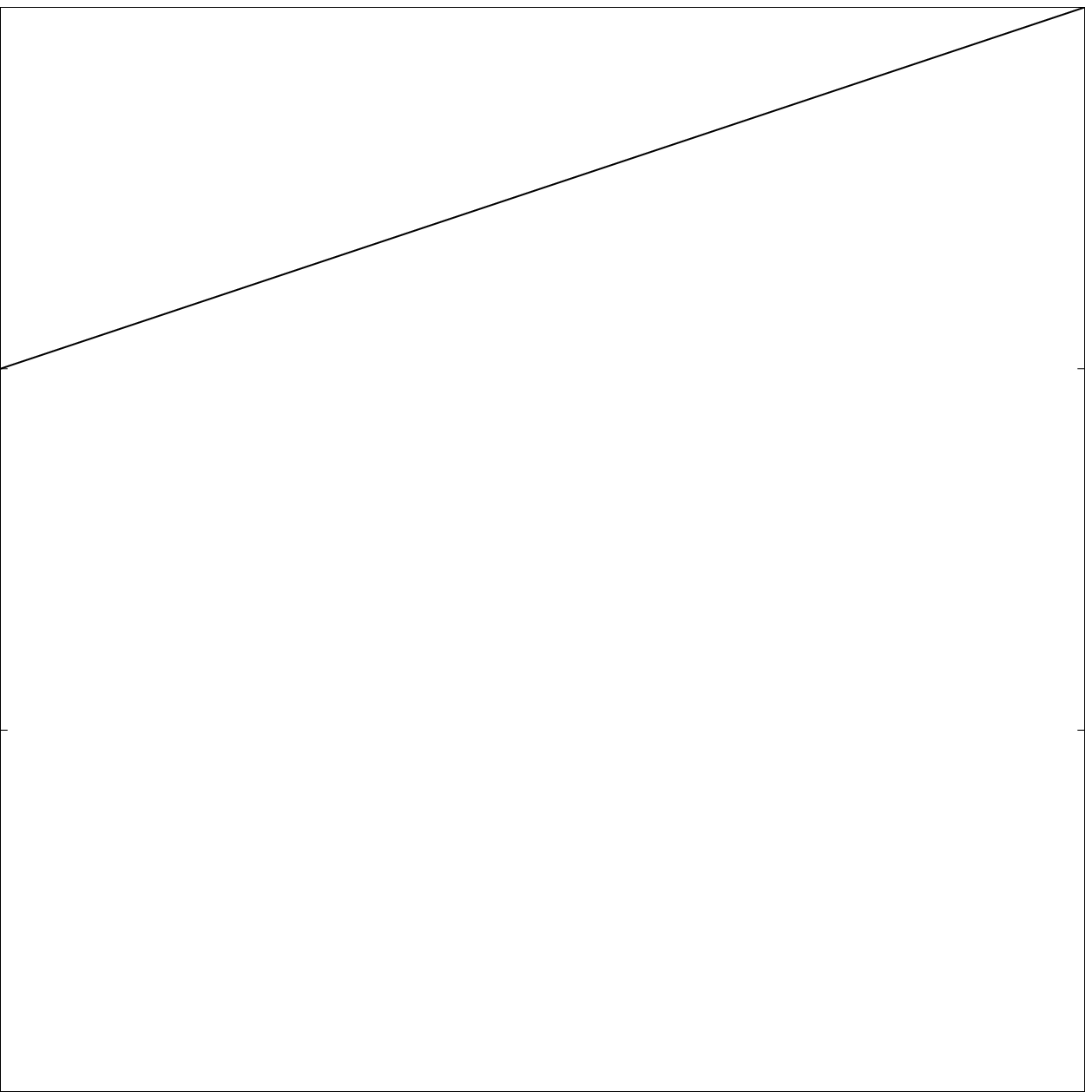}}
\put(-2,-2){\makebox(0,0)[tr]{$0$}}
\put(360,-2){\makebox(0,0)[t]{$1$}}
\put(-2,360){\makebox(0,0)[r]{$1$}}
\put(-2,120){\makebox(0,0)[r]{$\frac{1}{3}$}}
\put(-2,240){\makebox(0,0)[r]{$\frac{2}{3}$}}
\put(180,-62){\makebox(0,0)[t]{$\tau_{2}$}}
\end{picture}
\end{center}
\caption{The Cantor set $X$.}%
\label{Fig3}%
\end{figure}

\begin{exmp}
\label{ExaRep.5}\upshape Let $X$ be the middle-third Cantor set. The mapping
$\sigma\left(  x\right)  =3x\operatorname*{mod}1$ restricts to $X$ when $X$ is
embedded in the unit interval in the usual fashion; see Fig.\ \ref{Fig3}. The
two maps $\tau_{1}\left(  x\right)  =\frac{x}{3}$, $\tau_{2}=\frac{x+2}{3}$
satisfy $X=\tau_{1}\left(  X\right)  \cup\tau_{2}\left(  X\right)  $, and%
\begin{equation}
\sigma\circ\tau_{i}=\operatorname*{id}\nolimits_{X}. \label{eqRep.21}%
\end{equation}
The geometry of $X$ is illustrated in Fig.\ \ref{Fig3}. The Cantor measure
$\mu$ is determined uniquely by the two properties%
\begin{equation}
\mu\left(  X\right)  =1\text{\qquad and\qquad}\mu=\frac{1}{2}\left(  \mu
\circ\tau_{1}^{-1}+\mu\circ\tau_{2}^{-1}\right)  . \label{eqRep.22}%
\end{equation}
In fact, we have
\begin{equation}
\mu\circ\tau_{i}^{-1}=2\mu|_{\tau_{i}\left(  X\right)  }, \label{eqRep.23}%
\end{equation}
and (\ref{eqRep.10}) are clearly satisfied. The difference between the two
examples is that now $\mu$ is the Cantor measure, while in the previous
example it was the Lebesgue measure. The support of the Cantor measure is the
Cantor set $X$.

Using the corollary, we note that two isometries $S_{i}\colon L^{2}\left(
\mu\right)  \rightarrow L^{2}\left(  \mu\right)  $ are defined by the
formulas
\begin{equation}
S_{i}\varphi\left(  x\right)  =\varphi\left(  \sigma\left(  x\right)  \right)
\sqrt{2}\chi_{\tau_{i}\left(  X\right)  }\left(  x\right)  ,\label{eqRep.24}%
\end{equation}
where now $\sigma\left(  x\right)  =3x\operatorname{mod}1$, and%
\[
\tau_{1}\left(  X\right)  =X\cap\left[  0,\frac{1}{3}\right]  \text{\qquad
and\qquad}\tau_{2}\left(  X\right)  =X\cap\left[  \frac{2}{3},1\right]  .
\]
The formulae for the adjoint co-isometries $S_{i}^{\ast}\colon L^{2}\left(
\mu\right)  \rightarrow L^{2}\left(  \mu\right)  $ are
\begin{equation}
S_{1}^{\ast}\varphi\left(  x\right)  =\frac{1}{\sqrt{2}}\varphi\left(
\frac{x}{3}\right)  ,\qquad S_{2}^{\ast}\varphi\left(  x\right)  =\frac
{1}{\sqrt{2}}\varphi\left(  \frac{x+2}{3}\right)  ,\label{eqRep.25}%
\end{equation}
and it is immediate that the Cuntz relations (\ref{eqRep.12}) and
(\ref{eqRep.14}) are satisfied: by direct verification, or by an application
of Corollary \ref{CorRep.3}, we note that the isometries (\ref{eqRep.24}) form
a representation of $\mathcal{O}_{2}$ on the Hilbert space $L^{2}\left(
\mu\right)  $.
\end{exmp}

\section{\label{Fro}From representations to iterated function systems}

In Section \ref{Iso} we showed that every iterated function system (IFS), even
if not contractive, may be represented in the Hilbert space $\mathcal{H}$ of
equivalence classes $\operatorname*{class}\left(  \varphi,\mu\right)  $, where
$\varphi$ is a Borel function and $\mu$ is a positive Borel measure. The
equivalence relation $\sim$ which defines $\mathcal{H}$ is given by
(\ref{eqIso.1}). In Section \ref{Rep} we specialized this construction to the
case when measures $\mu$ may be found such that%
\begin{equation}
\mu\circ\tau_{i}^{-1}\absolutelycontinuous\mu, \label{eqFro.1}%
\end{equation}
where $\tau_{i}$, $i=1,\dots,N$, is the given IFS. In that case, we proved
that the resulting representation of $\mathcal{O}_{N}$ may be realized in
$L^{2}\left(  \mu\right)  $. For each such measure $\mu$ satisfying
(\ref{eqFro.1}), the representation of $\mathcal{O}_{N}$ on $L^{2}\left(
\mu\right)  $ is a sub-representation of the \textquotedblleft
global\textquotedblright\ representation on the Hilbert space $\mathcal{H}$
from Section \ref{Iso}.

In this section, the tables are turned: now the starting point is some
representation of $\mathcal{O}_{N}$, and we wish to reconstruct some IFS and
its realization in Hilbert space.

\begin{defn}
\label{DefFro.1}\upshape Let $N\in\mathbb{N}$, $N\geq2$, and let
$\mathcal{O}_{N}$ be the Cuntz algebra on $N$ generators, i.e., the $C^{\ast}%
$-algebra based on the relations (\ref{eqRep.12}) and (\ref{eqRep.14}). It is
known, see \cite{Cun77}, to be a simple $C^{\ast}$-algebra. If $\mathcal{K}$
is a Hilbert space, we say that $\mathcal{O}_{N}$ is represented on
$\mathcal{K}$ if there are isometries $S_{i}\colon\mathcal{K}\rightarrow
\mathcal{K}$ such that%
\begin{equation}
S_{i}^{\ast}S_{j}=\delta_{i,j}I_{\mathcal{K}}\text{\qquad and\qquad}\sum
_{i=1}^{N}S_{i}S_{i}^{\ast}=I_{\mathcal{K}}.\label{eqFro.2}%
\end{equation}
When $\mathcal{K}$ is given at the outset, we will denote all the
representations of $\mathcal{O}_{N}$ on $\mathcal{K}$ by $\operatorname*{Rep}%
\left(  \mathcal{O}_{N},\mathcal{K}\right)  $.

If $X$ is a compact Hausdorff space, then we denote by $\mathcal{H}\left(
X\right)  $ the Hilbert space of equivalence classes $\operatorname*{class}%
\left(  \varphi,\mu\right)  $ introduced in Section \ref{Iso}. If further
$\sigma,\tau_{i}\colon X\rightarrow X$ is a given IFS of continuous maps
satisfying the three conditions%
\begin{align}
&  \vphantom{\bigcup\nolimits_{i=1}^{N}}\sigma\circ\tau_{i}=\operatorname*{id}%
\nolimits_{X},\label{eqFro.3}\\
&  \vphantom{\bigcup\nolimits_{i=1}^{N}}\text{the maps }\tau_{i}\text{ are
distinct branches of the inverse for }\sigma\colon X\rightarrow X\text{,
and}\label{eqFro.3bis}\\
&  \bigcup\nolimits_{i=1}^{N}\tau_{i}\left(  X\right)  =X, \label{eqFro.4}%
\end{align}
then the isometries%
\begin{equation}
S_{i}\colon\left(  \varphi,\mu\right)  \longrightarrow\left(  \varphi
\circ\sigma,\mu\circ\tau_{i}^{-1}\right)  , \label{eqFro.5}%
\end{equation}
or%
\[
S_{i}\left(  \varphi\sqrt{d\mu}\right)  :=\varphi\circ\sigma\;\sqrt{d\mu
\circ\tau_{i}^{-1}},
\]
define an element in $\operatorname*{Rep}\left(  \mathcal{O}_{N}%
,\mathcal{H}\left(  X\right)  \right)  $, and we say that this is the
\emph{universal} representation of $\mathcal{O}_{N}$ built on the IFS $\left(
X,\tau_{i}\right)  $.

Our next result justifies this terminology.

If $\mathcal{K}$ and $\mathcal{H}$ are Hilbert spaces which both carry
representations of $\mathcal{O}_{N}$, we say that the representation on
$\mathcal{K}$ is a \emph{sub-representation} of that on $\mathcal{H}$ if there
is an isometry $W\colon\mathcal{K}\rightarrow\mathcal{H}$ such that%
\begin{equation}
WS_{i}^{\mathcal{K}}=S_{i}^{\mathcal{H}}W.\label{eqFro.6}%
\end{equation}
The isometry $W$ is said to \emph{intertwine} the two representations. (The
superscripts in the formula (\ref{eqFro.6}) indicate the Hilbert space on
which the isometries $S_{i}$ act.)
\end{defn}

\begin{thm}
\label{ThmFro.2}Let $N\in\mathbb{N}$, $N\geq2$, be given. Let $\left(
S_{i}\right)  $ be in $\operatorname*{Rep}\left(  \mathcal{O}_{N}%
,\mathcal{K}\right)  $ for some Hilbert space $\mathcal{K}$, and let $\left(
X,\sigma,\tau_{i}\right)  $ be an iterated function system on a compact metric
space $X$ which satisfies the conditions of Definition \textup{\ref{DefFro.1}%
}, and furthermore%
\begin{equation}
\operatorname*{diameter}\left(  \tau_{i_{1}}\circ\tau_{i_{2}}\circ\dots
\circ\tau_{i_{k}}\left(  X\right)  \right)  \underset{k\rightarrow\infty
}{\longrightarrow}0. \label{eqFro.7}%
\end{equation}
Then $\left(  S_{i},\mathcal{K}\right)  $ is a sub-representation of the
universal representation of $\mathcal{O}_{N}$ on $\mathcal{H}\left(  X\right)
$, i.e., there is an isometry $W\colon\mathcal{K}\rightarrow\mathcal{H}\left(
X\right)  $ such that%
\begin{equation}
WS_{i}^{\mathcal{K}}=S_{i}W, \label{eqFro.8}%
\end{equation}
where the isometries $S_{i}\colon\mathcal{H}\left(  X\right)  \rightarrow
\mathcal{H}\left(  X\right)  $ are defined as in \textup{(\ref{eqFro.5})}.
\end{thm}

\begin{pf}
Let $N$, $\mathcal{K}$, $X$, $\sigma$, $\tau_{i}$ be given as in the statement
of the theorem, and let $\left(  S_{i}\right)  \in\operatorname*{Rep}\left(
\mathcal{O}_{N},\mathcal{K}\right)  $. We shall omit the superscript when
referring to the isometries $S_{i}$ ($=S_{i}^{\mathcal{K}}$). It is easy to
see that for every $k\in\mathbb{N}$, the $N^{k}$ distinct projections
\begin{equation}
P\left(  i_{1},\dots,i_{k}\right)  :=S_{i_{1}}S_{i_{2}}\cdots S_{i_{k}%
}S_{i_{k}}^{\ast}\cdots S_{i_{2}}^{\ast}S_{i_{1}}^{\ast} \label{eqFro.9}%
\end{equation}
are mutually orthogonal, and that%
\begin{equation}
\sum_{i_{1},\dots,i_{k}}P\left(  i_{1},\dots,i_{k}\right)  =I_{\mathcal{K}}
\label{eqFro.10}%
\end{equation}
and%
\begin{equation}
\sum_{j=1}^{N}P\left(  i_{1},\dots,i_{k},j\right)  =P\left(  i_{1},\dots
,i_{k}\right)  . \label{eqFro.11}%
\end{equation}
Using an argument from \cite{Jor04a}, it follows that there is a unique
projection-valued measure $P\left(  \,\cdot\,\right)  $ defined on the Borel
sets in $X$ such that%
\begin{equation}
P\left(  \tau_{i_{1}}\circ\tau_{i_{2}}\circ\dots\circ\tau_{i_{k}}\left(
X\right)  \right)  =P\left(  i_{1},\dots,i_{k}\right)  , \label{eqFro.12}%
\end{equation}
where the notation is abused: $P$ is denoting both the expression in
(\ref{eqFro.9}) and the new measure. Specifically, $P\left(  \,\cdot\,\right)
$ satisfies the following additional five properties:

\begin{enumerate}
\item \label{ThmFro.1proof(1)}%
$\displaystyle\vphantom{\int_{X}\sum_{i=1}^{N}}I_{\mathcal{K}}=P\left(
X\right)  =\int_{X}P\left(  dx\right)  $;

\item \label{ThmFro.1proof(2)}%
$\displaystyle\vphantom{\int_{X}\sum_{i=1}^{N}}P\left(  E\right)  =P\left(
E\right)  ^{\ast}=P\left(  E\right)  ^{2}$,\quad$\displaystyle E\in
\mathcal{B}\left(  X\right)  $;

\item \label{ThmFro.1proof(3)}%
$\displaystyle\vphantom{\int_{X}\sum_{i=1}^{N}}P\left(  \,\cdot\,\right)  $ is
countably additive;

\item \label{ThmFro.1proof(4)}%
$\displaystyle\vphantom{\int_{X}\sum_{i=1}^{N}}P\left(  E\right)  P\left(
F\right)  =0$\quad if $\displaystyle E,F\in\mathcal{B}\left(  X\right)  $ and
$\displaystyle E\cap F=\varnothing$;

\item \label{ThmFro.1proof(5)}%
$\displaystyle\vphantom{\int_{X}\sum_{i=1}^{N}}\sum_{i=1}^{N}S_{i}P\left(
\tau_{i}^{-1}\left(  E\right)  \right)  S_{i}^{\ast}=P\left(  E\right)
$,\quad$\displaystyle E\in\mathcal{B}\left(  X\right)  $.
\end{enumerate}

As a result, we note that for every $f,g\in\mathcal{K}$, $\mu_{f,g}\left(
\,\cdot\,\right)  =\left\langle \,f\mid P\left(  \,\cdot\,\right)
g\,\right\rangle $ is a signed Borel measure on $X$. We shall use the
abbreviation $\mu_{f}\left(  E\right)  :=\mu_{f,f}\left(  E\right)
=\left\Vert P\left(  E\right)  f\right\Vert ^{2}$, $f\in\mathcal{K}$,
$E\in\mathcal{B}\left(  X\right)  $, and we note that $\mu_{f}$ is positive.
Moreover%
\begin{equation}
\left\vert \mu_{f,g}\left(  E\right)  \right\vert ^{2}\leq\mu_{f}\left(
E\right)  \mu_{g}\left(  E\right)  ,\qquad f,g\in\mathcal{K},\;E\in
\mathcal{B}\left(  X\right)  . \label{eqFro.13}%
\end{equation}

Formula (\ref{ThmFro.1proof(5)}) above specializes to the recursive identity%
\begin{equation}
\sum_{i=1}^{N}\mu_{S_{i}^{\ast}f}\circ\tau_{i}^{-1}=\mu_{f}.\label{eqFro.14}%
\end{equation}
Note that this is the covariance condition (\ref{eqInt.5}) from the
introduction. Substituting $S_{i}f$ for $f$, we get%
\begin{equation}
\mu_{f}\circ\tau_{i}^{-1}=\mu_{S_{i}f}.\label{eqFro.15}%
\end{equation}

We are now ready to define the operator $W\colon\mathcal{K}\rightarrow
\mathcal{H}\left(  X\right)  $ which intertwines the given representation of
$\mathcal{O}_{N}$ on $\mathcal{K}$ with the universal $\mathcal{O}_{N}%
$-representation acting on $\mathcal{H}\left(  X\right)  $:%
\[
W\colon\mathcal{K}\ni f\longmapsto\mu_{f}\longmapsto\left(  1,\mu_{f}\right)
\in\mathcal{H}\left(  X\right)  .
\]
It follows from the definitions that%
\begin{equation}
\left\Vert f\right\Vert _{\mathcal{K}}^{2}=\left\Vert \left(  1,\mu
_{f}\right)  \right\Vert _{\mathcal{H}\left(  X\right)  }^{2},
\label{eqFro.16}%
\end{equation}
where the norm $\left\Vert \,\cdot\,\right\Vert _{\mathcal{H}\left(  X\right)
}$ is defined in Section \ref{Iso}. Using the polarization identity%
\begin{equation}
\left\langle \,f\mid g\,\right\rangle =\frac{1}{4}\sum_{k=0}^{3}%
i^{k}\left\Vert i^{k}f+g\right\Vert ^{2}, \label{eqFro.17}%
\end{equation}
we further note that $W$ is an isometry from $\mathcal{K}$ to $\mathcal{H}%
\left(  X\right)  $. Moreover, formula (\ref{eqFro.15}) yields%
\begin{equation}
WS_{i}^{\mathcal{K}}f=\left(  1,\mu_{f}\circ\tau_{i}^{-1}\right)  =S_{i}Wf,
\label{eqFro.18}%
\end{equation}
where $S_{i}$ denotes the isometry in $\mathcal{H}\left(  X\right)  $ defined
in (\ref{eqIso.2}), i.e., the universal representation. \qed
\end{pf}

\begin{exmp}
\label{ExaFro.3}\upshape To see that there are examples which are covered by
Theorem \ref{ThmFro.2} but not by the more restrictive setting in Section
\ref{Rep}, we only need to specify an $\left(  S_{i}\right)  $ in, say,
$\operatorname*{Rep}\left(  \mathcal{O}_{2},\mathcal{K}\right)  $ and some
$e_{0}\in\mathcal{K}$ such that the measure $\mu:=\mu_{e_{0}}$ does not
satisfy the conditions (\ref{eqRep.10}). Here is a simple one: let
$\mathcal{K}:=L^{2}\left(  \mathbb{T}\right)  $, where $\mathbb{T}$ is the
one-torus $\mathbb{T}=\left\{  \,z\in\mathbb{C}\mid\left\vert z\right\vert
=1\,\right\}  $ equipped with Haar measure. Set%
\begin{equation}
S_{0}f\left(  z\right)  =f\left(  z^{2}\right)  ,\qquad S_{1}f\left(
z\right)  =zf\left(  z^{2}\right)  . \label{eqFro.19}%
\end{equation}
Then it is immediate that%
\begin{equation}
S_{i}^{\ast}S_{j}=\delta_{i,j}I\text{\qquad and\qquad}\sum S_{i}S_{i}^{\ast
}=I, \label{eqFro.20}%
\end{equation}
so (\ref{eqFro.19}) defines an element in $\operatorname*{Rep}\left(
\mathcal{O}_{2},L^{2}\left(  \mathbb{T}\right)  \right)  $. We take the
corresponding IFS to be the unit interval with the subdivision from Example
\ref{ExaRep.4}; see also Fig.\ \ref{Fig2}. Setting $e_{n}\left(  z\right)
=z^{n}$, we note that%
\[
S_{0}e_{n}=e_{2n},\qquad S_{1}e_{n}=e_{2n+1},\qquad n\in\mathbb{Z},
\]
so the two isometries permute the vectors in an orthonormal basis for
$\mathcal{K}=L^{2}\left(  \mathbb{T}\right)  $. If $f\in\mathcal{K}$, the
measures $\mu_{f}\left(  \,\cdot\,\right)  =\left\Vert P\left(  \,\cdot
\,\right)  f\right\Vert ^{2}$ are Borel measures on $\left[  0,1\right]  $.
Taking $f=e_{0}$, one easily checks that $\mu_{e_{0}}=\delta_{0}$, $\mu
_{e_{0}}\circ\tau_{0}^{-1}=\delta_{0}$, and $\mu_{e_{0}}\circ\tau_{1}%
^{-1}=\delta_{1/2}$ where $\delta_{0}$ and $\delta_{1/2}$ are the Dirac
measures at $x=0$ and $x=1/2$, respectively, and $\tau_{k}\left(  x\right)
:=\frac{x+k}{2}$, $k=0,1$. This makes it clear that (\ref{eqRep.10}) is not
satisfied. However, Theorem \ref{ThmFro.2} does apply to this example.
\end{exmp}

\begin{rem}
\label{RemFro.4}\upshape One might wonder how big a part of the universal
Hilbert space $\mathcal{H}\left(  X\right)  $ is needed for realizing the
representations of $\mathcal{O}_{N}$. The answer is that the most general
vectors $\varphi\sqrt{d\mu}$ are needed if we want to represent all the
elements in $\operatorname*{Rep}\left(  \mathcal{O}_{N},\mathcal{K}\right)  $.
To see this, note that if $\left(  S_{i}\right)  $ is in $\operatorname*{Rep}%
\left(  \mathcal{O}_{N},\mathcal{K}\right)  $ for some Hilbert space
$\mathcal{K}$, then the projections $S_{i_{1}}\cdots S_{i_{k}}S_{i_{k}}^{\ast
}\cdots S_{i_{1}}^{\ast}$, $k\in\mathbb{N}$, $i_{1},\dots,i_{k}\in
\mathbb{Z}_{N}$, generate a commutative $C^{\ast}$-algebra $\mathfrak{A}$ of
operators on $\mathcal{K}$. We showed in \cite{BJO04} that the corresponding
representations of $\mathfrak{A}$ include all possible spectral types. At the
same time, Nelson showed in \cite{Nel69} that the spectral representation,
including the multiplicity function, for abelian $C^{\ast}$-algebras
$\mathfrak{A}$ may be realized concretely in the Hilbert space $\mathcal{H}%
\left(  X\right)  $ where $X$ is the compact Gelfand space of $\mathfrak{A}$.
Hence all possible positive measures on $X$ will be needed in our
understanding of the representations of $\mathcal{O}_{N}$: the representations
of $\mathcal{O}_{N}$ may be realized acting on vectors $\varphi\sqrt{d\mu}%
\in\mathcal{H}\left(  X\right)  $, and all positive Borel measures on $X$ are
needed for this.
\end{rem}

\section{\label{Con}Conclusions}

This paper studies a general class of iterated function systems (IFS). No
contractivity assumptions are made, other than the existence of some compact
attractor. The possibility of escape to infinity is considered.

We are concerned with the realization of point transformations in Hilbert
space. Our present approach in fact is based directly on a certain
Hilbert-space construction, and on the theory of representations of the Cuntz
algebras $\mathcal{O}_{N}$, $N = 2, 3, \dots$.

While the more traditional approaches to IFS's start with some equilibrium
measure, ours doesn't. Rather, we construct a Hilbert space directly from a
given IFS, and our construction uses instead families of measures. Starting
with a fixed IFS, $\mathcal{S}_{N}$ with $N$ branches, we prove existence of
an associated representation of $\mathcal{O}_{N}$, and we show that the
representation is universal in a certain sense. Our framework includes as a
special case representations of $\mathcal{O}_{2}$ associated with quadrature
mirror filters and wavelets, and similarly, subband filters with $N$ frequency subbands.

We further prove a theorem about a direct correspondence between a given
system $\mathcal{S}_{N}$, and an associated sub-representation of the
universal representation of $\mathcal{O}_{N}$.%

\begin{ack}%
Our work is motivated by the lectures and papers by Arveson \cite{Arv03b} and
Tsirelson \cite{Tsi03} from the symposium \textquotedblleft Advances in
Quantum Dynamics\textquotedblright\ held at Mount Holyoke in 2002, sponsored
jointly by the AMS, IMS, and the SIAM. The paper \cite{BJO04} further offers
more results on the general theory of representations of the Cuntz algebras.
We thank Ola Bratteli for discussions. These papers make other interesting
applications of the Hilbert spaces $\mathcal{H}\left(  X\right)  $. They also
lists yet other applications of the idea underlying the $\mathcal{H}\left(
X\right)  $ construction, and the book \cite{Nel69} by Nelson gives a
beautiful presentation of the Spectral Representation Theorem, also couched in
terms of the same Hilbert spaces.

The Mount Holyoke symposium was supported by NSF grant DMS-9973450, and the
present research was supported by the NSF/FRG grant DMS-0139473.

We are grateful to Brian Treadway for corrections and suggestions, for
beautiful typesetting, and for graphics construction.
\end{ack}%


\end{document}